\documentclass[12pt]{amsart}
\usepackage{amsmath,amssymb,amscd,array, mathrsfs }

\usepackage{amsmath,amscd,amssymb,amsthm,array}
\usepackage{color}

\usepackage{amsmath,amssymb,mathrsfs,amsthm, tikz-cd,mathrsfs}
\usepackage{comment}
\def\url#1{\expandafter\s

\tring\csname #1\endcsname}

\def\mmat #1,#2,#3,#4,{\text{\small\arraycolsep=3pt $
\begin{pmatrix}#1&#2\\#3&#4\end{pmatrix}$}}

\usepackage{hyperref}
\usepackage{capt-of}

\usepackage{multirow}
\usepackage[all]{xy}

\usepackage{comments}
\newComments\SBe{Said}{blue}
\newComments\SBo{Sofiane}{blue}
\newComments\AM{Nacer}{blue}
\newComments\DL{DL}{red}
\newComments\QEh{QEh}{blue}

\def\mmat #1,#2,#3,#4,{\text{\small\arraycolsep=3pt $
\begin{pmatrix}#1&#2\\#3&#4\end{pmatrix}$}}

\usepackage{lscape}
\usepackage{tikz-cd}
\usepackage{enumerate}

\usepackage{DLdef1}
\usepackage{multicol}

\def\mmat #1,#2,#3,#4,{\text{\small\arraycolsep=3pt $
\begin{pmatrix}#1&#2\\#3&#4\end{pmatrix}$}}

\renewcommand {\ssbegin}[2][*]
 {\refstepcounter{subsection}%
\if#1*
\addcontentsline{toc}{subsection}{\thesubsection.\ #2}%
\else
\addcontentsline{toc}{subsection}{\thesubsection.\hskip 1pc #2. #1}%
\fi
 \def \secno {\gdef \secno {}{\ssecfont
\thesubsection.\hskip 2ex}%
 }%
 \begin{#2}}

\renewcommand {\sssbegin}[2][*]
 {\refstepcounter{subsubsection}
\if#1*
\addcontentsline{toc}{subsubsection}{\thesubsubsection.\hskip 1pc #2}%
\else
\addcontentsline{toc}{subsubsection}{\thesubsubsection.\hskip 1pc #2. #1}
\fi
 \def \secno {\gdef \secno {}{\ssecfont \thesubsubsection.\hskip 2ex}%
 }%
 \begin{#2}}

\renewcommand {\parbegin}[2][*]
 {\refstepcounter{paragraph}
\if#1*
\addcontentsline{toc}{paragraph}{\theparagraph.\hskip 1pc #2}%
\else
\addcontentsline{toc}{paragraph}{\theparagraph.\hskip 1pc #2. #1}
\fi
 \def \secno {\gdef \secno {}{\ssecfont \theparagraph.\hskip 2ex}%
 }%
 \begin{#2}}

\rmnameii{etr}{etr}
\rmnameii{evv}{ev}
\DeclareMathOperator{\h}{\mathcal{H}}

\DeclareMathOperator{\K}{\mathbb{K}}

\newcommand{\Ll}{{\mathrm{L}}}
\newcommand{\Rr}{{\mathrm{R}}}

    \newcommand{\Om}{\Omega}

\newcommand{\al}{\alpha}

\newcommand{\la}{\lambda}

\newcommand{\de}{\delta}
\newcommand{\om}{\omega}

\def\br{[\;,\;]}

\newcommand{\g}{\mathfrak{g}}
\begin{document}

\title[Structure of Flat Quadratic Quasi-Frobenius Flat Lie Superalgebras]{Structure of Flat Quadratic Quasi-Frobenius Lie Superalgebras via Double Extensions}

\author{Sofiane Bouarroudj}

\address{Division of Science and Mathematics, New York University Abu Dhabi, P.O. Box 129188, Abu Dhabi, United Arab Emirates.}
\email{sofiane.bouarroudj@nyu.edu}

\author{Hamza El Ouali}
\address{Universit\'e Cadi-Ayyad,
	Facult\'e des sciences et techniques,
	B.P. 549, Marrakech, Maroc.}
\email{eloualihamza11@gmail.com}

\thanks{S.B. was supported by the grant NYUAD-065}

\keywords{Quasi-Frobenius Lie (super)algebras, the Levi-Civita product, symplectic product, quadratic Lie (super)algebras, double extension.}

 \subjclass[2020]{17B05;17A60; 17A70}

\begin{abstract} 
A \emph{flat quadratic quasi-Frobenius} Lie superalgebra is a quadratic Lie superalgebra equipped with an additional symplectic structure that is flat with respect to the natural symplectic product. In this paper, we introduce the notion of a  \emph{flat quadratic double extension} of a flat quadratic quasi-Frobenius Lie superalgebra, in the cases where both the symplectic structure and the quadratic structure are either even or odd. We show that, over an algebraically closed field, any such Lie superalgebra can be constructed through a sequence of flat quadratic double extensions starting from the trivial algebra $\{0\}$. Moreover, when the quadratic and symplectic structures have different parity, we introduce the notion of a \emph{planar double extension}, which constitutes the main novelty of this paper. In this case, we prove that such Lie superalgebras have  total dimension $4n$.

Finally, we  classify flat quadratic quasi-Frobenius non-abelian Lie superalgebras of dimension at most four and present explicit examples in dimensions six and eight.
\end{abstract}


\maketitle

\thispagestyle{empty}
\setcounter{tocdepth}{2}
\section{Introduction}

\subsection{Double extensions of quadratic Lie (super)algebras}
The notion of double extension of Lie algebras equipped with a  nondegenerate, invariant (or associative) {\it symmetric}  bilinear form was introduced by Medina and Revoy in \cite{MR}, although earlier instances already appeared in \cite{K1}. Such Lie algebras are commonly referred to as {\it quadratic} Lie algebras. 

The concept of double extension was developed to characterize solvable Lie algebras and, in particular, was applied to the study of nilpotent Lie algebras. More generally, it was shown that any quadratic Lie algebra with non-trivial center can be described inductively from another quadratic Lie algebra of co-dimension $2$ (see also \cite{FS}). Other interesting examples include Manin triples, Manin–Olshansky triples, and affine Kac–Moody Lie algebras with Cartan matrix. Each of these structures can be realized as a double extension.

The first superization of the results of Medina and Revoy is due to Benamor and Benayadi in \cite{BEBE}. The authors generalized this notion to even quadratic Lie superalgebras and characterized solvable Lie superalgebras using this concept. Later, Bajo-Benayadi-Bordemann introduced the notion of generalized double extension for quadratic Lie superalgebras \cite{BBB}, providing an inductive description and classification of such superalgebras. In \cite{ABB}, the authors introduced the notion of  generalized double extension for odd-quadratic Lie superalgebras and provided an inductive description of solvable odd-quadratic Lie superalgebras, as well as of odd-quadratic Lie superalgebras whose even part is a reductive Lie algebra. This study turns out to be a nontrivial superization, since many classical tools—such as the Levi decomposition and Lie’s theorem—do not hold in the super setting. 


\subsection{Double extensions of quasi-Frobenius Lie (super)algebras}
 Following \cite{BE, BEh,  BM, BR}, a Lie  superalgebra is called orthosymplectic (resp. periplectic) quasi-Frobenius Lie superalgebra if it is endowed with a nondegenerate closed orthosymplectic (resp. periplectic) \emph{antisymmetric} bilinear form, see Section \ref{QFslags}
 
 Medina and Revoy introduced the symplectic double extension of Lie algebras in \cite{MR1}. Roughly speaking, the idea is to enlarge a quasi-Frobenius Lie algebra by adding a symplectic plane, thereby producing a new quasi-Frobenius Lie algebra.  This construction serves as the antisymmetric analogue of the double extension of quadratic Lie algebras formulated by Medina and Revoy in \cite{MR}. The symplectic double extension was developed to characterize and classify nilpotent symplectic Lie groups through their Lie algebras. Recall that a Lie group admits an invariant symplectic structure if it carries a left-invariant and closed $2$-form of maximal rank. However, the original computation of the symplectic double extension overlooked a necessary cocycle condition. Medina and Dardi\'e later corrected and generalized the construction in \cite{MD}, where the proper formulation was established. Their work showed that every nilpotent quasi-Frobenius Lie algebra can be obtained by successive symplectic double extensions starting from the trivial Lie algebra. 
 
 In \cite{DM}, the authors introduced the notion of Kähler double extension for Kähler Lie algebras, which realizes a Kähler Lie algebra as a Kähler reduction of another one. They showed that every Kähler Lie algebra can be obtained by a finite sequence of such extensions starting either from the trivial algebra $\{0\}$ or from a flat Kähler Lie algebra. Moreover, they proved that any completely solvable and unimodular Kähler Lie algebra is necessarily abelian. 
 
In \cite{BM}, the authors generalized the double extension method to the super setting by studying Lie superalgebras endowed with nondegenerate closed orthosymplectic or periplectic forms. They proved that every quasi-Frobenius Lie superalgebra satisfying certain structural conditions can be obtained as a double extension of a smaller quasi-Frobenius Lie superalgebra.
Further in \cite{BEh, F}, the authors developed the process of symplectic double extensions for Lie (super)algebras with  degenerate centers.
\subsection{The natural symplectic product}
Let $(M, \Omega)$ be a symplectic manifold, that is, a manifold endowed with a non-degenerate closed 2-form. A symplectic connection is a connection $\nabla$ such that $\nabla \Omega = 0$. In contrast with pseudo-Riemannian manifolds, which admit a unique torsion-free connection preserving the metric (the Levi-Civita connection), a symplectic manifold possesses many torsion-free symplectic connections (see \cite{BCGRS}).

In \cite{BE}, it was shown that a quasi-Frobenius (symplectic) Lie superalgebra $(\mathfrak{g}, \omega)$ always admits a symplectic product, although this product is not unique. A bilinear product $\star$ is called symplectic if
$$
[u, v] = u \star v - (-1)^{|u||v|} v \star u, \quad \text{and} \quad
\omega(u \star v, w) + (-1)^{|u||v|} \omega(v, u \star w) = 0,
$$
for all homogeneous $u,v,w \in \mathfrak{g}$. In contrast, for pseudo-Euclidean non-associative superalgebras, the Levi-Civita product exists and is unique.  However, in \cite{BE}, the authors introduced a natural symplectic product on Lie superalgebras that depends uniquely on the Lie bracket and the symplectic form. This product is defined by
$$
\omega(u \star v, w) = \frac{1}{3} (\omega([u,v], w) + (-1)^{|v||w|}\omega([u,w], v)), \quad \text{for all homogeneous } u,v,w \in \mathfrak{g}.
$$
The curvature of $(\mathfrak{g}, \omega)$ at the identity is defined by (where $\Ll_u$ is the left multiplication):
$$
K(u,v) := \Ll_{[u,v]} - [\Ll_u, \Ll_v], \quad \text{for all homogeneous } u,v \in \mathfrak{g}.
$$

A quasi-Frobenius Lie superalgebra $(\mathfrak{g}, \omega)$ is called flat if its curvature $K$ vanishes identically. In this case, flatness is equivalent to $(\mathfrak{g}, \star)$ being a left-symmetric algebra. 

Flat quasi-Frobenius Lie algebras were studied in \cite{BELL}, while \cite{BE} investigated flat quasi-Frobenius Lie superalgebras equipped with a flat natural symplectic product. In \cite{BE, BELL}, the notion of flat double extension was introduced, and it was shown that every flat quasi-Frobenius Lie (super)algebra can be obtained through a sequence of flat double extensions starting from $\{0\}$; moreover, any such Lie (super)algebra is necessarily nilpotent.
 \subsection{Double extensions of quadratic quasi-Frobenius Lie (super)algebras} 
A quasi-Frobenius (or symplectic) quadratic Lie superalgebra 
is a Lie superalgebra endowed with both a quadratic structure 
and a quasi-Frobenius (or symplectic) structure. In this paper, we will refer to them as QQF-superalgebras. 

In \cite{BBM}, Bajo-Benayadi-Medina studied  QQF-algebras and provided inductive descriptions of such Lie algebras using the framework of double extensions.
In \cite{BaBe}, the authors introduced various notions of double extensions for even quadratic quasi-Frobenius Lie superalgebras—when both structures are even—as well as for Manin triples, thereby establishing inductive construction methods.


Finally, in \cite{BBNS2}, the authors conducted an in-depth study of a sub-class of QQF-superalgebras. They proved, in particular, that all such superalgebras are nilpotent. As a consequence, they provided classifications in low dimensions and identified the double extensions that preserve symplectic structures. Using both elementary odd double extensions and generalized double extensions of QQF-superalgebras, they obtained an inductive description of QQF-superalgebras of filiform type.

\subsection{Outline of the results}
In the present work, we focus on flat QQF-superalgebras, which are flat with respect to the natural symplectic product. Note that flat quadratic pseudo-Riemannian Lie algebras were previously studied in \cite{BRM} in the Lorentzian case, and subsequently in \cite{BL}, where the notion of double extension was introduced. In particular, it was shown that every flat quadratic Lorentzian Lie algebra can be obtained from a Euclidean abelian algebra by means of double extensions.

In this paper, we prove that a flat quasi-Frobenius Lie superalgebra $(\mathfrak{g}, \omega)$ admits a quadratic structure if and only if there exists an invertible $\omega$-antisymmetric endomorphism $\rho$ such that, for all $u \in \mathfrak{g}$, the endomorphism $\rho \circ \mathrm{ad}_u$ is $\omega$-symmetric, see Prop.  \ref{prop:quadratic-structure}.  Furthermore, we show that in this case the total dimension of the superalgebra is necessarily even, see Prop. \ref{dim-even}. 

When $\rho$ is even, we adapt the notion of flat double extension  to the setting of flat QQF-superalgebras and show that every flat QQF-superalgebra can be obtained as a flat double extension of another flat QQF-superalgebra of smaller dimension, see Theorem \ref{quadratic-ex1} and Theorem \ref{quadratic-ex2}. In particular, any flat QQF-superalgebra arises from the trivial algebra $\{0\}$ through a finite sequence of flat double extensions, see Corollary \ref{suite-ex} and Corollary \ref{lcro}.

When $\rho$ is odd, the classical flat double extension cannot be adapted. To overcome this difficulty, we introduce the notion of flat planar double extension for flat QQF-superalgebras, see Sec. \ref{oddrho}. We prove that every flat QQF-superalgebra with an odd $\rho$ can be obtained by a flat {\it planar} double extension of a flat QQF-superalgebra of smaller dimension, see Theorem \ref{thm-rho-od1} and Theorem \ref{thm-rho-od2}. Furthermore, we show that in this case the dimension of the superalgebra is necessarily of the form $4n$, where $n \in \mathbb{N}$, and that every flat QQF-superalgebra with odd $\rho$ can be obtained from the trivial algebra $\{0\}$ by a finite sequence of flat planar double extensions (see Corollary~\ref{cosuite}).

Finally, in Section~\ref{examples-low}, we classify flat QQF-superalgebras of dimension 
$4$ and present several explicit examples in dimensions $6$ and 
$8$.

\section{Backgrounds}\label{back}
 Let $\mathbb K$ be an algebraically closed field of characteristic zero, although most of the constructions remain valid over any field of characteristic greater than 
3. We write $\mathbb{Z}_2=\{\bar 0, \bar 1\}$ for the group of integers modulo 2. For a more comprehensive treatment of linear algebra over superspaces, the reader is referred to \cite{L, S}.

Let $V=V_{\bar 0}\oplus V_{\bar 1}$ be a superspace defined over  ${\mathbb K}$. The parity of a homogeneous element $v\in V_{\bar{i}}$ is denoted by $|v|:=\bar{i}$. The element $v$ is called \textit{even} if $v\in V_{\bar 0}$ and \textit{odd} if $v\in V_{\bar 1}$. 

Throughout the text, all elements in  linear expression are assumed to be homogeneous, unless otherwise stated. 

Let \(V\) be a vector superspace. For homogeneous elements \(x,y\in V\) we denote by \(x^*,y^*\in V^*\) their duals . For every \(x^*\otimes y^*\in V^*\otimes V^*\) and  \(u\otimes v\in V\otimes V\), we define the pairing
\[
\langle x^*\otimes y^*, u \otimes v\rangle :=(-1)^{|y||u|}x^*(u)\,y^*(v).
\]
The wedge product (graded antisymmetrisation) is given by
\[
x^*\wedge y^*=x^*\otimes y^* - (-1)^{|x||y|}\,y^*\otimes x^*.
\]
The symmetric product (graded symmetrisation) is given by
\[
x^*\odot y^*=x^*\otimes y^* + (-1)^{|x||y|}\,y^*\otimes x^*.
\]

A linear map $\varphi:V\rightarrow W$ between two superspaces is called \textit{even} if $\varphi(V_{\bar{i}})\subset W_{\bar{i}}$ and \textit{odd} if $\varphi(V_{\bar{i}})\subset W_{\overline{i}+\overline{1}}$.

\subsection{Bilinear forms on a superspace}
Let $V$ and $W$ be two superspaces. Let ${\cal B}\in \text{Bil}(V,W)$ be a homogeneous bilinear form. Recall that the Gram matrix $B=(B_{ij})$ associated to $\cal B$ is given by the following formula:
\begin{equation*}\label{martBil}
B_{ij}=(-1)^{|{\cal B}| |v_i|}{\cal B}(v_{i}, w_{j})\text{~~for the basis vectors $v_{i}\in V$ and $w_{j}\in W$ .}
\end{equation*}
This definition allows us to identify a~bilinear form of $\mathrm{Bil}(V, W)$ with an element of $\mathrm{Hom}(V, W^*)$. 
Consider the \textit{upsetting} of bilinear forms
$u:\text{Bil} (V, W)\rightarrow \text{Bil}(W, V)$ given by the formula 
\begin{equation*}\label{susyB}
u({\cal B})(w, v)=(-1)^{|v| | w|}{\cal B}(v,w)\text{~~for any $v \in V$ and $w\in W$.}
\end{equation*}
In terms of the Gram matrix $B$ of ${\cal B}$, we have 
\begin{equation*}
u(B)=
\left ( \begin{array}{cc}  R^{t} & (-1)^{|B|}T^{t} \\ (-1)^{|B|}S^{t} & -U^{t} \end{array} \right ),
\text{ for $B=\left (  \begin{array}{cc}R & S \\ T & U \end{array} \right )$}.
\end{equation*}
From now on, we will assume that $V=W$. Following \cite{L}, we say that: 
\begin{itemize} \item the form
$\cal B$ is  \textit{symmetric} if  and only if $u(B)=B$;

\item the form $\cal B$ is \textit{anti-symmetric} if and only if $u(B)=-B$.

\end{itemize}
Given a linear map $f: V\rightarrow V$. The adjoint of $f$ with respect to ${\cal B}$ is the linear map $f^*$ defined as follows:
\[
{\cal B}(f(v),w)= (-1)^{|f||v|}{\cal B}(v,f^*(w)).
\]
We say that $f$ is ${\cal B}$-{\it symmetric} (resp. ${\cal B}$-{\it antisymmetric}) if $f^*=f$ (resp. $f^*=-f$).

Throughout the text, we denote by $\omega$ an {\it antisymmetric} bilinear form, meaning that for all 
$u,v \in V$, 
\[
\omega(u,v)=-(-1)^{|u||v|} \omega(v,u).
\]
A non-degenerate, anti-symmetric even form is called {\it orthosymplectic}; an odd one is called {\it periplectic}. Hereafter, we assume that such forms are homogeneous. Following \cite{BM}, if $\omega$ is periplectic, then 
$\dim(V_{\bar 0}) = \dim(V_{\bar 1})$, 
whereas if $\omega$ is orthosymplectic, then $\dim(V_{\bar 0})$ is even.

Throughout the text, we denote by ${\mathscr B}$ a  {\it symmetric} bilinear form, meaning that for all 
$u,v \in V$, 
\[
{\mathscr B}(u,v)=(-1)^{|u||v|} {\mathscr B}(v,u). 
\]
Hereafter, we assume that such forms are homogeneous. 
\subsection{Lie superalgebras with symmetric and anti-symmetric bilinear forms}\label{QFslags}
Let $\fg$ be a (finite-dimensional) Lie superalgebra over $\mathbb K$. Following \cite{BE, BM}, a Lie  superalgebra $\g$ is called orthosymplectic (resp. periplectic) quasi-Frobenius Lie superalgebra if it is endowed with an orthosymplectic (resp. periplectic) form  $\omega$ that is closed. That is, 
\begin{equation}
(-1)^{|u||w|}\omega( [u, v], w ) + (-1)^{|v||u|} \omega( [v, w], u )+  (-1)^{|w||v|}\omega([w, u], v ) = 0, \text{ for all } u,v,w \in \g.
\label{cyclic}
\end{equation}
We denote such a Lie superalgebra  by $(\g, \omega)$. The condition \eqref{cyclic} means that $\omega$ is a 2-cocycle on $\fg$ with values in $\K$. If $\omega$ is a coboundary, then $(\g, \omega)$ is called a {\it Frobenius} Lie superalgebra, see \cite{O}. Examples in low dimensions can be found in \cite{BM,BE}, while \cite{BR} investigates such structures for all real Lie superalgebras of total dimension~4.

Let \((\g, \om)\) be a quasi-Frobenius Lie superalgebra. It was proved in \cite[Prop. 3.18]{BE} that there exists a unique {\it natural symplectic product} that depends on the product  \(\br\) and the form \(\om\). This product is given by: 
\begin{equation}
\om( u \star v, w ) = \frac{1}{3} \left( \om( [u, v], w ) +(-1)^{|v||w|} \om( [u, w], v ) \right), \quad \text{ for all \(u, v, w \in \g\).}
\label{cyclic13}\end{equation} 
The product $\star$ in Eq. (\ref{cyclic13}) was first introduced in \cite{BELL} in the context of Lie algebras. It verifies two important properties: (where $\mathrm{L}_u^\star$ denotes the left multiplication by $u$):
\[
[u,v]=u\star v- (-1)^{|u||v|} v \star u, \quad  \text{ and } \quad  (\mathrm{L}_u^\star)^*=- \mathrm{L}^\star_u, \quad  \text{ for all $u,v\in \g$. }
\]
This product differs from the so-called Levi--Civita product; see \cite{BE} for further details.

The product $\star$ in Eq. (\ref{cyclic13}) is called {\it flat} if
\[
\Ll^\star_{[u,v]} = [\Ll^\star_u, \Ll^\star_v], \qquad \text{for all $u,v \in \g$.}
\]
A quasi-Frobenius Lie superalgebra $(\mathfrak{g}, \omega)$ is called \emph{flat quasi-Frobenius} if the \emph{natural quasi-Frobenius product} associated with $(\g,\omega)$ is flat. The structure of these Lie superalgebras has recently been studied in \cite{BE}. For instance, these Lie superalgebras are necessarily nilpotent in the  case where $\omega$ is homogeneous; see \cite{BE}.

A Lie superalgebra $\g$ is called quadratic\footnote{This structure is referred to as a nis in \cite{BKLS}.} if it is equipped with a nondegenerate, invariant, symmetric bilinear form $\mathscr B$. That is
\[
{\mathscr B}([u,v],w)={\mathscr B}(u, [v,w]) \quad  \text{ for all } u,v,w\in \g.
\]
We denote such a Lie superalgebra  by $(\g, {\mathscr B})$. 

Finite-dimensional Lie superalgebras over $\mathbb C$ with symmetrizable Cartan matrix are quadratic, see \cite{K1}. Some other simple Lie superalgebras which are not subquotients of Lie superalgebras with Cartan matrix are quadratic; for instance, $\mathfrak{h}^{(1)}(0|m)$ and  $\mathfrak{psq}(n)$ are quadratic, see \cite{BKLS} for more details.


Now, a Lie superalgebra $\mathfrak{g}$ is called \emph{quadratic quasi-Frobenius} if it is equipped with a quadratic structure ${\mathscr B}$ and a quasi-Frobenius  structure $\omega$. For the sake of conciseness, we will refer to them as $\mathrm{QQF}$-superalgebras. We denote such a structure by $(\mathfrak{g}, \omega, {\mathscr B})$.


 

A $\mathrm{QQF}$-superalgebra $(\fg, \omega, \mathscr B)$ is called flat if the natural symplectic product \eqref{cyclic13} associated with $(\g,\omega)$ is flat. 
 \section{Quadratic vs. Quasi-Frobenius Structures on Lie Superalgebras}\label{quad-symp}
The purpose of this section is to demonstrate that quadratic and quasi-Frobenius structures are related via an invertible derivation. In the case of Lie algebras, Jacobson \cite{J} proved that any finite-dimensional Lie algebra over a field of characteristic zero that admits an invertible derivation is necessarily nilpotent. On the other hand, we have already shown in \cite{BE} that any flat quasi-Frobenius Lie superalgebra is nilpotent.

Fix a quadratic structure. The following proposition describes the conditions under which a quasi-Frobenius structure exists.

\ssbegin{Proposition}\label{CQ}
Let $(\g, {\mathscr B})$ be a quadratic Lie superalgebra. 
A quasi-Frobenius structure $\omega$ can be defined on $\g$ 
if and only if $\g$ admits an invertible ${\mathscr B}$-antisymmetric derivation $\delta$. 
Moreover, $|\delta|=| {\mathscr B}|+|\omega|$. 
\end{Proposition}
\begin{proof}
Since $ {\mathscr B}$ is non-degenerate, there exists an invertible homogeneous endomorphism $\delta$  such that $
\omega(u,v)= {\mathscr B}(\delta(u),v),$  for all $u,v\in\g.$ Since  $ {\mathscr B}$ is symmetric and $\omega$ is anti-symmetric, it follows that $\delta$ is $ {\mathscr B}$-anti-symmetric. For all $u,v,w\in \g$, we have 
\begin{align*}
&(-1)^{|u||w|}\omega(u,[v,w]) + (-1)^{|u||v|}\omega(v,[w,u]) + (-1)^{|v||w|}\omega(w,[u,v]) \\
&= (-1)^{|u||w|} {\mathscr B}(\delta(u),[v,w]) + (-1)^{|u||v|} {\mathscr B}(\delta(v),[w,u]) + (-1)^{|v||w|} {\mathscr B}(\delta(w),[u,v]) \\
&= (-1)^{|u||w|} {\mathscr B}([\delta(u),v],w) + (-1)^{|u||w|+|\delta||u|} {\mathscr B}([u,\delta(v)],w) - (-1)^{|u||w|} {\mathscr B}(\delta([u,v]),w).
\end{align*}
This shows that $\omega$ is closed if and only if 
\[
\delta([u,v]) = [\delta(u),v] + (-1)^{|\delta||u|}[u,\delta(v)].
\]
From the construction of $\delta$, we deduce that   $|\delta|=|{\mathscr B}|+|\omega|$.
\end{proof}

\ssbegin{Remark}
It should be noticed that under the assumptions of Prop. \ref{CQ}, the derivation $\delta$ of $\mathfrak{g}$ is already assumed to be  ${\mathscr B}$-antisymmetric.  Moreover, if $\delta$ is even, then it is also $\omega$-antisymmetric, whereas if $\delta$ is odd, it becomes $\omega$-symmetric.
Indeed, for all $u, v \in \mathfrak{g}$,
\[
\omega(\delta u, v) =  {\mathscr B}(\delta^2(u), v)
= -(-1)^{|\de|(|\delta|+|u|)}  {\mathscr B}(\delta u, \delta v)
= -(-1)^{|\de|+|\delta||u|} \omega(u, \delta v).
\]
It is natural to conjecture that every quasi-Frobenius Lie superalgebra $(\g,\omega)$ with an invertible, $\omega$-antisymmetric derivation carries a quadratic structure. Nevertheless, we will demonstrate that this does not hold in general. For instance, the Lie superalgebra $(\mathfrak{g}^3, \omega)$ defined  \textup{(}in the basis $\{x_1, x_2\mid  y_1, y_2\}$\textup{)} by:
$$[x_1, y_1]=y_2, \quad \om= x_1^*\wedge y_2^*+x_2^* \wedge y_1^*,$$  is quasi-Frobenius (see, \cite{BE}) but does not admit a quadratic structure. Choose specifically the linear map $\delta$ of $\mathfrak{g}$ defined by
\[
\delta= x_1\otimes x_1^* +2 x_2 \otimes x_2^*-2y_1\otimes y_1^*- y_2 \otimes y_2^*.
\]
It is easy to check that $\delta$ is an $\omega$-antisymmetric invertible derivation.  

For the sake of contradiction, suppose that $\mathfrak{g}^3$ admits a quadratic structure $ {\mathscr B}$. Then we would have
\[
1=\omega(x_1, y_2)  = \omega(x_1, [x_1, y_1]) =  {\mathscr B}(\delta(x_1), [x_1, y_1]) =   {\mathscr B}(x_1, [x_1, y_1]) =   {\mathscr B}([x_1, x_1], y_1) = 0,
\]
which is a contradiction.
\end{Remark}

\ssbegin{Example}\label{ex1} Following \cite{BE}, the Lie superalgebra $(\g^2, \omega)$ defined as:  \textup{(}in the basis $\{x_1, x_2\mid  y_1, y_2\}$\textup{)}
$$
	[x_1,y_1]=y_2,\quad [y_1, y_1]= x_2,\; \mbox{ and } \; \om=2x_1^*\wedge x_2^*-y_1^*\wedge y_2^*,
	$$
is a flat QQF-superalgebra. The even symmetric bilinear form $ {\mathscr B}$ on $\g^2$ is defined by 
\[
 {\mathscr B}=-x_1^* \odot x_2^*-y_1^*\odot y_2^*,
\]
while the even invertible derivation  $\delta$ on $\g^2$ is defined by
\[
\delta= -2x_1 \otimes x_1^*+2x_2 \otimes x_2^*+y_1 \otimes y_1^*-y_2\otimes y_2^*.
\] 
\end{Example}

\ssbegin{Example}\label{ex2} Following \cite{BE}, the Lie superalgebra $(\g^4, \omega)$ defined as: 
\textup{(}in the basis $\{x_1, x_2\mid  y_1, y_2\}$\textup{)}
$$[y_1, y_1]=x_1,\quad [y_1, y_2]=x_2,  \; \text{ and } \; \om= -2x_1^*\wedge y_2^*+x_2^* \wedge y_1^*,$$
is a flat QQF-superalgebra. The odd symmetric bilinear form $ {\mathscr B}$ on $\g^4$ is defined by 
\[
 {\mathscr B}=x_1^*\odot y_2^*+x_2^*\odot y_1^*,
\]
while the even invertible derivation $\delta$ on $\g^4$ is defined by
\[
\delta= -2 x_1 \otimes x_1^*+ x_2\otimes x_2^*-y_1 \otimes y_1^*+2y_2\otimes y_2^*.
\]
\end{Example}

We now present a method for constructing flat QQF-superalgebras. As we shall see, this class is quite rich and abundant, which further motivates and justifies a systematic study of these structures.

Let $(\mathfrak{h}, \omega_\mathfrak{h}, {\mathscr B}_\mathfrak{h})$ be a flat QQF-superalgebra, and let $(\mathcal{A} , \cdot, \Omega)$ be  an associative supercommutative   superalgebra, endowed with an even nondegenerate invariant symmetric bilinear form $\Omega$, i.e.,
\[
\Omega(a\cdot b, c) = \Omega(a, b\cdot c), \quad  \text{ for all } \quad a,b,c \in \mathcal{A}.
\]
On the superspace $\mathfrak{g} = \g_{\bar 0}\oplus \g_{\bar 1}$, where $\g_{\bar 0}=\mathfrak{h}_{\bar 0}\otimes \mathcal{A}_{\bar 0} \oplus \mathfrak{h}_{\bar 1}\otimes \mathcal{A}_{\bar 1}$  and $\g_{\bar 1}=\mathfrak{h}_{\bar 0}\otimes \mathcal{A}_{\bar 1} \oplus \mathfrak{h}_{\bar 1}\otimes \mathcal{A}_{\bar 0}$, we define  the following structures:
\begin{equation}\label{Om-g}
\begin{array}{ccl}
[u\otimes a, v \otimes b] &:=& (-1)^{|a||v|}[u,v]_\mathfrak{h} \otimes (a \cdot b), \\[2mm]
\omega(u \otimes a, v \otimes b) &:=& (-1)^{|a||v|}\omega_\mathfrak{h}(u, v) \, \Omega(a,b), \;\,\\[2mm]
{\mathscr B}(u \otimes a, v \otimes b) & := & (-1)^{|a||v|} {\mathscr B}_\mathfrak{h}(u, v)\, \Omega(a,b),
\end{array}
\end{equation}
for all $u,v \in \mathfrak{h}$ and $a,b \in \mathcal{A}$.   
\ssbegin{Proposition} The superalgebra defined by \textup{Eq.} \eqref{Om-g} is a flat QQF-superalgebra.
\end{Proposition}
\begin{proof}
It is easy to see that $(\mathfrak{g},\om, {\mathscr B})$ is a  QQF-superalgebra. We now show that the natural symplectic product associated with  $(\mathfrak{g}, \omega)$ is flat.  

Let $(u \otimes a) \star (v \otimes b)$ (resp. $u \star_\mathfrak{h} v$) denote the natural symplectic product of $(\mathfrak{g}, \omega)$ (resp. $(\mathfrak{h}, \omega_\mathfrak{h})$). For any $u,v,w \in \mathfrak{h}$ and any $a,b,c \in \mathcal{A}$, we have \small
\begin{align*}
&\omega((u \otimes a) \star (v \otimes b), w \otimes c)\\  
&\quad = \frac{1}{3} \Big( \omega([u \otimes a, v \otimes b], w \otimes c) 
+ (-1)^{|v\otimes b||w \otimes c|} \omega([u \otimes a, w \otimes c], v \otimes b)\Big)\\
&\quad = \frac{1}{3} \Big( (-1)^{|a||v|}\omega([u , v ]_\mathfrak{h}\otimes (a\cdot b), w \otimes c) 
+ (-1)^{|v\otimes b||w \otimes c|+|a||w|} \omega([u , w ]_\mathfrak{h}\otimes(a\cdot c), v \otimes b)\Big) \\ 
&\quad = \frac{1}{3} \Big((-1)^{|a||v|+(|a|+|b|)|w|} \omega([u , v ]_\mathfrak{h},w )\Om(a\cdot b, c) 
+ (-1)^{X} \omega([u , w ]_\mathfrak{h}, v )\Om(a\cdot c, b)\Big) \\ 
&\quad = \frac{1}{3} \Big( (-1)^{|a||v|+(|a|+|b|)|w|}\omega([u , v ]_\mathfrak{h},w )\Om(a\cdot b, c) 
+ (-1)^{X+|b||c|} \omega([u , w ]_\mathfrak{h}, v )\Om(a\cdot b, c)\Big) \\ 
&\quad = \frac{1}{3}(-1)^{|a||v|+(|a|+|b|)|w|} \Big( \omega([u , v ]_\mathfrak{h},w )
+ (-1)^{|v||w |} \omega([u , w ]_\mathfrak{h}, v )\Big) \Om(a\cdot b, c) \\   
&\quad = (-1)^{|a||v|+(|a|+|b|)|w|} \omega(u\star_\mathfrak{h} v,w )
 \Om(a\cdot b, c) \\
&\quad= (-1)^{|a||v|}\om((u \star_\mathfrak{h} v) \otimes (a \cdot b), w \otimes c),
\end{align*}
\normalsize
where we have put $X=|v\otimes b||w \otimes c|+|a||w|+ (|a|+|c|)|v|$.
We conclude that 
\[
(u \otimes a) \star (v \otimes b) =(-1)^{|a||v|} (u \star_\mathfrak{h} v) \otimes (a \cdot b).
\]
Since the natural symplectic product $\star_\mathfrak{h}$ in $(\mathfrak{h}, \omega_\mathfrak{h})$ is left-symmetric and $\mathcal{A}$ is associative and supercommutative, it follows that the natural symplectic product in $(\mathfrak{g}, \omega)$ is also left-symmetric. Thus, $(\g, \om, {\mathscr B})$ is a flat QQF-superalgebra. 
\end{proof}



Contrariwise to Prop. \ref{CQ}, let us now fix a quasi-Frobenius structure. The proposition below describes the conditions under which a quadratic structure exists.

\ssbegin{Proposition}\label{prop:quadratic-structure}
Let $(\mathfrak{g}, \omega)$ be a quasi-Frobenius Lie superalgebra \textup{(}not necessarily flat\textup{)}. A quadratic structure $ {\mathscr B}$ can be defined on $\mathfrak{g}$ if and only if there exists an 
 invertible endomorphism $\rho$ of $\mathfrak{g}$ such that $\rho$ is $\om$-antisymmetric and the operator $\rho \circ \operatorname{ad}_u$ is $\omega$-symmetric, for all $u \in \mathfrak{g}$. Moreover,  $|\rho|= | {\mathscr B}|+|\omega|$. 
\end{Proposition}
\begin{proof}
Let $(\mathfrak{g}, \omega)$ be a quasi-Frobenius Lie superalgebra.  
Assume that $ {\mathscr B}$ is a non-degenerate symmetric bilinear form on $\mathfrak{g}$.  
There exists an invertible homogeneous endomorphism $\rho$ of $(\mathfrak{g}, \omega)$ such that
$
 {\mathscr B}(-,-) = \omega(\rho(-), -)$. Since $\omega$ is anti-symmetric and  $ {\mathscr B}$ is symmetric , it follows that $\rho$ is $\om$-anti-symmetric. For all $u,v,w \in \mathfrak{g}$, we have
\[
 {\mathscr B}([u,v],w) + (-1)^{|u||v|}  {\mathscr B}(v,[u,w])
= \omega\big(\rho \circ \operatorname{ad}_u(v), w\big)
- (-1)^{|v|(|\rho|+|u|)} \omega\big(v, \rho \circ \operatorname{ad}_u(w)\big).
\]
Hence, $ {\mathscr B}$ is invariant if and only if 
the operator $\rho \circ \operatorname{ad}_u$ is $\omega$-symmetric, for any $u \in \mathfrak{g}$.
\end{proof}

\ssbegin{Remark}
\label{rem:rho-delta}
(i) Sometimes, we denote a flat QQF-superalgebra by 
\((\mathfrak{g}, \omega, \rho)\) instead of \((\mathfrak{g}, \omega,  {\mathscr B})\).  
Note that the endomorphism $\rho$ is also $ {\mathscr B}$-anti-symmetric.

(ii) If we fix a quadratic structure $ {\mathscr B}$ on $\mathfrak{g}$,  
Prop.~\ref{CQ} implies that  there exists an invertible anti-symmetric 
endomorphism $\delta$ of $(\mathfrak{g}, {\mathscr B})$ defined by
$
\omega(-,-) =  {\mathscr B}(\delta(-), -).$ We can also express the natural symplectic product associated with $(\mathfrak{g}, \omega)$ 
in terms of $\delta$ and the Lie bracket. Indeed, it is given by
\begin{equation}\label{NPL}
u \star v = \tfrac{1}{3}[u,v] + (-1)^{|u||\delta|}\, \delta^{-1}([u, \delta(v)]).
\end{equation}
The flatness of $(\mathfrak{g}, \omega)$ is equivalent to the fact that 
the product $\star$ is left-symmetric. 

(iii) Note that $\delta=\rho^{-1}$
where $\rho$ is the operator introduced in Prop.~\ref{prop:quadratic-structure}.  
In this paper, we choose to fix the flat quasi-Frobenius Lie superalgebra 
$(\mathfrak{g}, \omega)$ and express the quadratic structure by means of $\rho$.
Equivalently, one could fix the quadratic Lie superalgebra $(\mathfrak{g},  {\mathscr B})$ 
and describe the quasi-Frobenius structure using $\delta$ and the left-symmetric product 
defined in \eqref{NPL}.
\end{Remark}

\ssbegin{Proposition}\label{roh-inv}
Let $(\g,\omega,\rho)$ be a flat  QQF-superalgebra.  Then
\begin{enumerate}
    \item[$(i)$] $Z(\g)=[\g,\g]^{\perp_\omega}=[\g,\g]^{\perp_\rho}$,
    \item[$(ii)$] $\rho\bigl(Z(\g))=Z(\g)$.
\end{enumerate}
\end{Proposition}

\begin{proof}
From \cite[Prop. 3.3]{BE}, we have 
$
Z(\g)\subseteq[\g,\g]^{\perp_\omega}.
$
Since ${\mathscr B}(-,-)=\omega(\rho(-),-)$ defines a  quadratic form, it is easy to see that
$
[\g,\g]^{\perp_\rho}=Z(\g).
$
Moreover, since
\(\dim([\g,\g]^{\perp_\rho})=\dim([\g,\g]^{\perp_\omega})\), 
it follows that 
\[
Z(\g)=[\g,\g]^{\perp_\omega}=[\g,\g]^{\perp_\rho}.
\]
Let us prove Part (ii). If $a\in Z(\g)$, then
$
0=\omega(a,[u,v])={\mathscr B}    (\rho(a),[u,v]),
$
for all $u,v\in\g$. This shows that $\rho(a)\in Z(\g)$. 
\end{proof}

\ssbegin{Remark}
In \cite[Prop. 7.5]{BE}  we established the existence of exactly four non-abelian flat quasi-Frobenius Lie superalgebras of total dimension four. We have showed that $(\g^2, \omega)$ 
and $(\g^4, \omega)$ do admit a quadratic structure (see Examples \ref{ex1}  and \ref{ex2}). However, the other two Lie superalgebras do not admit such a structure. This is because $[\g, \g]^\perp \neq Z(\g)$, and the statement follows from Prop. \ref{roh-inv}.
\end{Remark}

\ssbegin{Proposition}\label{dim-even}
Let $(\g, \om, \rho)$ be a QQF-superalgebra.  Then the total dimension of $\g$ is even.\end{Proposition}
\begin{proof}
If $\om$ is odd, then from \cite[Lemma. 3.3]{BM} directly implies that $\dim(\g)$ is even. Assume now that $\om$ is even and the bilinear form
$
{\mathscr B}(-,-)=\om(\rho(-),-)$ is also even. Since $\om$ is anti-symmetric, its restriction to $\g_{\bar{0}}$ is also antisymmetric and nondegenerate. Therefore, the space $\g_{\bar{0}}$ must have even dimension.

Moreover, since $\mathscr B$ is symmetric, we have $
{\mathscr B}(u,v)=-\mathscr{B}(v,u),$ for all $u,v\in\g_{\bar{1}},$
and it is nondegenerate on $\g_{\bar{1}}$. Hence, $\dim(\g_{\bar{1}})$ is even. Therefore
$
\dim(\g)=\dim(\g_{\bar{0}})+\dim(\g_{\bar{1}})
$
is even.

Finally, if $\omega$ is even and $\mathscr B$ is odd, then by applying to $\mathscr B$ the same argument as in \cite[Lemma. 3.3]{BM}, we conclude that $\dim (\g)$ is even.
\end{proof}
\section{Flat double extensions}
The notion of double extension for flat quasi-Frobenius Lie superalgebras was introduced in \cite{BE}. For convenience, we recall this construction below, as it will be essential for formulating the corresponding double extension procedure for flat QQF–superalgebras. It is worth mentioning that the construction of double extensions 
for quadratic Lie superalgebras has been extensively studied by several authors 
(see, for instance, \cite{ABB}, \cite{BEBE} and references therein).

\ssbegin{Theorem}[See \cite{BE}]\label{double-ex1}
Let $(\mathfrak{b}, \omega_{\mathfrak{b}})$ be a flat orthosymplectic quasi-Frobenius Lie superalgebra, and let $\star_{\mathfrak{b}}$ be  the natural symplectic  product associated with it. 
Let $V = \mathbb{K}d$ be a one-dimensional vector space and $V^* = \mathbb{K}e$ be its dual. Assume that there exist an even linear map 
$\xi : \mathfrak{b} \to \mathfrak{b}$, 
and an element $b_{0} \in \mathfrak{b}_{\bar 0}$, satisfying the system \begin{equation}\label{eq:claim2}
\begin{array}{lcllcl}
\xi ([u,v]_\fb) & = &  \mathrm{L}_u^\star (\xi(v)) - (-1)^{|u||v|} \mathrm{L}_v^\star (\xi(u)), & 
\xi^* \circ \xi & = &  \tfrac{1}{3}(\Rr^\star_{b_0}+ (\Rr^\star_{b_0})^*), \\[2mm]
D(u\star v) & = & D(u)\star v+u\star D(v)-\xi(u)\star v , & 
[\xi, \xi^*] & = & \xi^2 - \tfrac{1}{3} \Rr^\star_{b_0},\;
b_0\in \ker{D},
\end{array}
\end{equation}
for all $u,v\in \mathfrak{b}$, where $D=\xi^*-\xi$. Define 
$
\g := \mathbb{K}d \oplus \mathfrak{b} \oplus \mathbb{K}e,
$
equipped with the Lie bracket and even bilinear form $\omega$ as follows:
\begin{equation}\label{Liebrackets1}
[d,u] = (\xi^* - 2\xi)(u) + \omega_{\mathfrak{b}}(b_0, u)e, \quad\text{ and }\quad
[u,v] = [u,v]_{\mathfrak{b}} + \omega_{\mathfrak{b}}\big((\xi + \xi^*)(u), v\big)e,
\end{equation}
for all $u,v \in \mathfrak{b}$ and
$$
\omega|_{\mathfrak{b} \times \mathfrak{b}} = \omega_{\mathfrak{b}},
\qquad \omega(e,d) = -\omega(d,e) = 1,
\qquad \omega(d,\mathfrak{b}) = \omega(e,\mathfrak{b}) = 0.
$$
Then $(\g, \omega)$ is a flat orthosymplectic quasi-Frobenius Lie superalgebra.
\end{Theorem}
Following \cite{BE}, we call  $(\g, \omega)$ the even flat double extension of $(\mathfrak{b}, \omega_{\mathfrak{b}})$ by means of $(\xi, b_0).$

\ssbegin{Theorem}[See \cite{BE}]\label{double-ex2}
Let $(\mathfrak{b}, \omega_{\mathfrak{b}})$ be a   flat orthosymplectic quasi-Frobenius Lie superalgebra, and let $\star_{\mathfrak{b}}$ denote the natural symplectic  product associated with it.
Let $V = \mathbb{K}d$ be a purely odd one-dimensional
vector superspace and $V^* = \mathbb{K}e$ be its dual. Assume that there exist an odd homogeneous linear map 
$\xi : \mathfrak{b} \to \mathfrak{b}$, 
and an even element $b_0 \in \mathfrak{b}_{\bar 0}$, satisfying the system \begin{equation}\label{eq:claim4}
\begin{array}{rclrcl}
\xi([u,v]_\fb) & = &  \mathrm{L}_u^\star (\xi(v)) - (-1)^{|u||v|} \mathrm{L}_v^\star (\xi(u)),&\left[\xi , \xi^*\right]
& = &  \Rr^\star_{b_0}-3\xi^2,\\[2mm]
 D(u\star v)
& = & D(u)\star v+(-1)^{|u|}u\star D(v)+\xi(u)\star v,& 
\Ll^\star_{b_0} & = & -(\xi+\xi^*)^2 ,\\[2mm]
(2\xi+\xi^*)(b_0)&=&0,\end{array}
\end{equation}
for all $u, v \in \mathfrak{b}$.   Define the vector superspace
$
\g := \mathbb{K}d \oplus \mathfrak{b} \oplus \mathbb{K}e,
$
equipped with the Lie bracket and even bilinear form $\omega$ as follows: 
\begin{equation}\label{Liebrackets2}
\begin{array}{lcllcl}
 [d,u] & = & - (-1)^{|u|}(\xi^* + 2\xi)(u) + \omega_{\mathfrak{b}}(b_0, u)e, & [d,d]& = & 2b_0,\\[2mm] 
[u,v] &=& [u,v]_{\mathfrak{b}} + \big((-1)^{|v|}\omega_{\mathfrak{b}}(\xi(u), v) +(-1)^{|u|}\omega_{\mathfrak{b}}(\xi^*(u),  v)\big)e,\end{array}
\end{equation}
for all $u,v \in \mathfrak{b}$ and
$$
\omega|_{\mathfrak{b} \times \mathfrak{b}} = \omega_{\mathfrak{b}},
\qquad \omega(e,d) = \omega(d,e) = 1,
\qquad \omega(d,\mathfrak{b}) = \omega(e,\mathfrak{b}) = 0.
$$
Then, $(\g, \omega)$ is a flat orthosymplectic quasi-Frobenius Lie superalgebra.

\end{Theorem}
Following \cite{BE}, we call $(\g, \omega)$ the odd flat double extension of $(\mathfrak{b}, \omega_{\mathfrak{b}})$ by means of $(\xi, b_0)$.

\ssbegin{Theorem}[See \cite{BE}]\label{double-ex3} Let $(\mathfrak{b}, \omega_{\mathfrak{b}})$ be a flat periplectic  quasi-Frobenius Lie superalgebra, and let $\star_{\mathfrak{b}}$ denote the natural symplectic  product associated with it.
Let $V = \mathbb{K}d$ be a one-dimensional vector space and  and $\Pi(V^*) = (\mathbb{K} e)_{\bar{1}}$, where $V^*$ is its dual. Assume that there exist an even linear map 
$\xi : \mathfrak{b} \to \mathfrak{b}$, 
and an even element $b_0 \in \mathfrak{b}_{\bar 0}$, satisfying the system \begin{equation}\label{eq:claim5}
\begin{array}{lcllcl}
\xi ([u,v]_\fb) & = &  \mathrm{L}_u^\star (\xi(v)) - (-1)^{|u||v|} \mathrm{L}_v^\star (\xi(u)), & 
\xi^* \circ \xi & = &  \tfrac{1}{3}(\Rr^\star_{b_0}+ (\Rr^\star_{b_0})^*), \\[2mm]
D(u\star v) & = & D(u)\star v+u\star D(v)-\xi(u)\star v , & 
[\xi, \xi^*] & = & \xi^2 - \tfrac{1}{3} \Rr^\star_{b_0},\;
b_0\in \ker{D},
\end{array}
\end{equation}
for all $u,v\in \mathfrak{b}$, where $D=\xi^*-\xi$. Define the vector superspace
$
\g :=V \oplus \mathfrak{b} \oplus \Pi(V^*),
$
equipped with the Lie bracket and odd bilinear form $\omega$ as follows:
\begin{equation}\label{Liebrackets3}
[d,u] = (\xi^* - 2\xi)(u) + \omega_{\mathfrak{b}}(b_0, u)e, \quad\text{ and }\quad
[u,v] = [u,v]_{\mathfrak{b}} + \omega_{\mathfrak{b}}\big((\xi + \xi^*)(u), v\big)e,
\end{equation}
for all $u,v \in \mathfrak{b}$, and
$$
\omega|_{\mathfrak{b} \times \mathfrak{b}} = \omega_{\mathfrak{b}},
\qquad \omega(e,d) = -\omega(d,e) = 1,
\qquad \omega(d,\mathfrak{b}) = \omega(e,\mathfrak{b}) = 0.
$$
Then, $(\g, \omega)$ is a flat periplectic  quasi-Frobenius Lie superalgebra.
\end{Theorem}
Following \cite{BE}, we call  $(\g, \omega)$ the even flat double extension of $(\mathfrak{b}, \omega_{\mathfrak{b}})$ by means of $(\xi, b_0)$.

\ssbegin{Theorem}[See \cite{BE}]\label{double-ex4}
Let $(\mathfrak{b}, \omega_{\mathfrak{b}})$ be a   flat periplectic quasi-Frobenius Lie superalgebra, and let $\star_{\mathfrak{b}}$ denote the natural symplectic  product associated with it. Let $V=V_{\bar{1}} = (\mathbb{K} d)_{\bar{1}}$ be a purely odd one-dimensional vector superspace  and $\Pi(V^*) = \mathbb{K} e$ where $V^*$ is its dual. Assume that there exist an odd homogeneous linear map 
$\xi : \mathfrak{b} \to \mathfrak{b}$, 
and an even element $b_0 \in \mathfrak{b}_{\bar 0}$, satisfying the system \begin{equation}\label{eq:claim6}
\begin{array}{rclrcl}
\xi([u,v]_\fb)&=&  \mathrm{L}_u^\star (\xi(v)) - (-1)^{|u||v|} \mathrm{L}_v^\star( \xi(u)),& [\xi ,\xi^*]  &=&  3\xi^2 + \Rr^\star_{b_0},\\[2mm] D(u\star v)&=& D(u)\star v+(-1)^{|u|}u\star D(v)-
(-1)^{|u|} \xi(u) \star v, &
\Ll^\star_{b_0}  &=& (\xi-\xi^*)^2,  \\[2mm] 
2\xi(b_0)&=&\xi^*(b_0),\end{array}
\end{equation}
for all $u, v \in \mathfrak{b}$. Define the vector superspace
$
\g := \mathbb{K}d \oplus \mathfrak{b} \oplus \mathbb{K}e,
$
equipped with the Lie bracket and odd bilinear form $\omega$ given as follows: 
\begin{equation}\label{Liebrackets4}
\begin{array}{lcllcl}
[d,u] & = & (-1)^{|u|}(\xi^* -2\xi)(u) + \omega_{\mathfrak{b}}(b_0, u)e, & [d,d] & = &-2b_0,\\ [2mm]
[u,v] & = &  [u,v]_{\mathfrak{b}} + \big((-1)^{|v|}\omega_{\mathfrak{b}}(\xi(u), v) +(-1)^{|u|}\omega_{\mathfrak{b}}(\xi^*(u),  v)\big)e,\end{array}
\end{equation}
for all $u,v \in \mathfrak{b}$ and
$$
\omega|_{\mathfrak{b} \times \mathfrak{b}} = \omega_{\mathfrak{b}},
\qquad \omega(e,d) = -\omega(d,e) = 1,
\qquad \omega(d,\mathfrak{b}) = \omega(e,\mathfrak{b}) = 0.
$$
Then $(\g, \omega)$ is a   flat periplectic  quasi-Frobenius Lie superalgebra.
\end{Theorem}
Following \cite{BE}, we call  $(\g, \omega)$ the odd flat double extension of $(\mathfrak{b}, \omega_{\mathfrak{b}})$ by means of $(\xi, b_0)$.

\ssbegin{Remark}
The terminology reflects the parity of the map $\xi$: the double extension is said to be even if $\xi$ is even, and odd if $\xi$ is odd. Moreover, since $\omega_\fb$ is homogeneous, no ambiguity arises.
\end{Remark}

\subsection{The case where \texorpdfstring{$|\rho|=\bar 0$}{|rho|=0}.}
We are now prepared to construct the notion of double extension for flat QQF–superalgebras in the case where $\rho$ is even, i.e., when the quadratic and quasi-Frobenius structures share the same parity.
\sssbegin{Theorem}\label{THQ1}
Let $(\mathfrak{b}, \omega_{\mathfrak{b}}, \rho_{\mathfrak{b}})$ be a flat  QQF-superalgebra, where $\rho_{\mathfrak{b}}$ is even. Let  
$(\mathfrak{g}, \omega)$ be the even flat double extension of  
$(\mathfrak{b}, \omega_{\mathfrak{b}})$ by means of $(\xi, b_0)$.  
Suppose that there exist $a_0 \in \mathfrak{b}_{\bar{0}}$ and $\lambda \in \mathbb{K}\backslash \{0\}$ such that 
\begin{equation}\label{DEX1}
\rho_{\mathfrak{b}} \circ ( 2\xi-\xi^{*} )+\lambda(\xi+\xi^*) = \Rr_{a_0}^{\star_\mathfrak{b}}+ (\Rr_{a_0}^{\star_\mathfrak{b}})^*, 
\qquad \text{and} \qquad (2\xi^* - \xi)(a_0) = \lambda\, b_0.
\end{equation}
For any $t \in \mathbb{K}$, we define the even linear map $\rho$ of $\mathfrak{g}$ by \textup{(}for all $u \in \mathfrak{b}$\textup{)}:
\[
\rho(e) = \lambda e, \qquad 
\rho(u) = \rho_{\mathfrak{b}}(u) + \omega_{\mathfrak{b}}(a_0, u)\, e, \qquad 
\rho(d) = \left \{ 
\begin{array}{ll}
t e + a_0 - \lambda d, & \mathrm{if }\; |\omega|=\bar 0,\\[2mm] a_0 - \lambda d, & \mathrm{if } \; |\omega|=\bar 1.
\end{array}
\right. 
\]
Then $(\mathfrak{g}, \omega, \rho)$ is a flat QQF-superalgebra, where $\rho$ is even.
\end{Theorem}
The flat QQF-superalgebra $(\g, \omega, \rho)$ is called the even flat double extension of $(\mathfrak{b}, \omega_{\mathfrak{b}}, \rho_\mathfrak{b})$ by means of $(\xi, b_0, a_0, \lambda)$. 
\begin{proof} According to Prop. ~\ref{prop:quadratic-structure}, it suffices to show that  $\rho$ is an invertible, $\omega$-antisymmetric endomorphism and that $\rho \circ \mathrm{ad}_a$ is $\omega$-symmetric, for any $a \in \mathfrak{g}$.  

Assume that $\omega$ is orthosymplectic. The proof in the periplectic case is completely analogous. Since $\rho_{\mathfrak{b}}$ is an even invertible $\omega_{\mathfrak{b}}$-antisymmetric endomorphism of 
$\mathfrak{b}$ and $\lambda \not =0$, it follows that 
$\rho$ is also an even invertible endomorphism of $\mathfrak{g}$  and also it is $\omega$-antisymmetric. The condition that $\rho \circ \mathrm{ad}_d$ is $\omega$-symmetric  
is equivalent to
\[
(2\xi^* - \xi)(a_0) = \lambda b_0, 
\qquad \text{and} \qquad 
\rho_{\mathfrak{b}} \circ (2\xi - \xi^{*})
\text{ is $\omega_{\mathfrak{b}}$-symmetric.}
\]
Similarly, the condition that $\rho \circ \mathrm{ad}_u$ is $\om$-symmetric for every $u \in \mathfrak{b}$ 
is equivalent to
\begin{equation}
\label{thm-rho}
\rho_{\mathfrak{b}} \circ ( \xi^{*}-2\xi ) 
= -(\Rr_{a_0}^{\star_\mathfrak{b}}+ (\Rr_{a_0}^{\star_\mathfrak{b}})^*) + \lambda(\xi + \xi^{*}).
\end{equation}
Note that the condition (\ref{thm-rho}) implies that $\rho_{\mathfrak{b}} \circ (2\xi - \xi^{*})$ is $\omega_{\mathfrak{b}}$-symmetric.  
Therefore, if the relations \eqref{DEX1} are satisfied, then 
$\rho \circ \mathrm{ad}_a$ is $\omega$-symmetric, for all $a \in \mathfrak{g}$.  It follows that the bilinear form $ {\mathscr B}$ defined by
$
 {\mathscr B}(-,-) = \omega(\rho(-), -),
$
is quadratic on 
$\mathfrak{g}$. Moreover, since both $\omega$ and $\rho$ are even, the quadratic form $ {\mathscr B}$ is even as well.
\end{proof}
We now adapt this method to the odd flat double extensions of flat quasi-Frobenius Lie superalgebras endowed with a quadratic structure.

\sssbegin{Theorem}\label{THQ2}
Let $(\mathfrak{b}, \omega_{\mathfrak{b}}, \rho_{\mathfrak{b}})$ be a flat QQF-superalgebra, where $\rho_{\mathfrak{b}}$ is even. Let  
$(\mathfrak{g}, \omega)$ be the odd flat double extension of  
$(\mathfrak{b},  \omega_{\mathfrak{b}})$ by means of $(\xi, b_0)$. Suppose that there exist $a_1 \in \mathfrak{b}_{\bar{1}}$ and $\lambda \in \mathbb{K}\backslash \{0\}$ such that 
\begin{itemize}
\item The case where $|\omega|=\bar 0$:
\begin{equation}\label{DEX2}
\begin{array}{rcl}
\rho_{\mathfrak{b}} \circ (2\xi + \xi^*)(u)  +\lambda (\xi - \xi^*)(u) & = & (-1)^{|u|}(\Rr_{a_1}^{\star_\mathfrak{b}}+(\Rr_{a_1}^{\star_\mathfrak{b}})^*)(u), \text{ and } \\[2mm] 
(\xi + 2\xi^{*})(a_1) + 2\,\rho_{\mathfrak b}(b_0) & = & -\lambda\, b_0.
\end{array}
\end{equation}
\item The case where $|\omega|=\bar 1$:
\begin{equation}\label{DEX3}
\begin{array}{rcl}
\rho_{\mathfrak{b}} \circ (2\xi - \xi^*)(u)  +\lambda (\xi + \xi^*)(u) & = & (-1)^{|u|}(\Rr_{a_1}^{\star_\mathfrak{b}}+(\Rr_{a_1}^{\star_\mathfrak{b}})^*)(u), 
\text{ and } \\[2mm]
(\xi-2\xi^*)(a_1)+2\rho_{\mathfrak{b}}(b_0) & = & - \lambda \, b_0.
\end{array}
\end{equation}
\end{itemize}
For any $t \in \mathbb{K}$, define the even linear endomorphism $\rho$ of $\mathfrak{g}$ by \textup{(}for all $u \in \mathfrak{b}$\textup{)}: 
\[
\rho(e) = \lambda e, \qquad 
\rho(u) = 
\rho_{\mathfrak{b}}(u)+\om_\mathfrak{b}(a_1, u) e,  \qquad 
\rho(d) = \left
\{
\begin{array}{ll}
t e -a_1- \lambda d & \text{if $|\omega|=\bar 0,$}\\[2mm]
a_1 -\lambda d & \text{if $|\omega|=\bar 1$.}
\end{array}
\right.
 \] 
Then $(\mathfrak{g}, \omega, \rho)$ is a flat  QQF-superalgebra, where $\rho$ is even.
\end{Theorem}

The flat QQF-superalgebra $(\g, \omega, \rho)$ is called the odd  flat double extension of $(\mathfrak{b}, \omega_{\mathfrak{b}}, \rho_\mathfrak{b})$ by means of $(\xi, b_0, a_1, \lambda)$.

\begin{proof}
The proof is similar to that of Theorem \ref{THQ1}.
\end{proof}


We now prove the converses of Theorems \ref{THQ1} and \ref{THQ2} when $|\omega|=\bar 0$.

\sssbegin{Theorem}\label{quadratic-ex1}
Let $(\g, \om, \rho)$ be a flat QQF-superalgebra such that $|\rho|=|\omega|=\bar 0$. Then $(\g, \om, \rho)$  is either:
\begin{itemize}
    \item[$(i)$] an even flat  double extension of a flat QQF-superalgebra $(\mathfrak{b}, \om_\mathfrak{b}, \rho_\mathfrak{b})$ by means of $(\xi, b_0, a_0, \lambda)$, or
    \item[$(ii)$] an odd flat  double extension of a flat QQF-superalgebra $(\mathfrak{b}, \om_\mathfrak{b}, \rho_\mathfrak{b})$ by means of $(\xi,\allowbreak b_0,\allowbreak a_1,\lambda)$.
\end{itemize}
\end{Theorem}

\begin{proof}
By \cite[Theorem 6.1]{BE}, we  have $Z(\g)\neq 0$.  Prop.~\ref{roh-inv} implies that 
$
\rho\bigl(Z(\g))
   =Z(\g).
$
Since $\K$ is algebraically closed, there exists 
$e\in Z(\g)$ such that $\rho(e)=\lambda e$ with $\lambda\neq 0$.

If $e$ is even, put   $I:=\mathbb{K}e$. This space is totally isotropic as 
$
\omega(e,e)=0.
$ By \cite[Theorem 6.7]{BE}, we deduce that  
$(\g,\omega)$ is an even flat double extension of  
$(\mathfrak{b},\omega_{\mathfrak{b}})$ 
by means of $(\xi,b_0)$. Let us write
$
\g=\K e\oplus\mathfrak{b}\oplus\K d,
$
where the Lie brackets are given by
\[
[d,u]
   =(\xi^{*}-2\xi)(u)+\omega_{\mathfrak{b}}(b_{0},u)\,e,
   \qquad
[u,v]=[u,v]_{\mathfrak{b}}+\om((\xi+\xi^{*})(u),v)\,e, \quad \text{for all $u,v\in\mathfrak{b}$.}
\] 
Since $\rho(I)=I$ and $\rho$ is $\omega$-antisymmetric, we have 
$
\rho(I^{\perp_\om})\subseteq I^{\perp_\om},$
 i.e.
$\rho(I\oplus\mathfrak{b})\subseteq I\oplus\mathfrak{b}.
$
It follows that $\rho$ can be written as
\[
\rho(e)=\lambda e,\qquad
\rho(u)=\rho_{\mathfrak{b}}(u)+\psi(u)e,\qquad
\rho(d)=t\,e+c_{0}+s\,d, \quad \text{for all $u\in\mathfrak{b}$,}
\]
 where  
$\rho_{\mathfrak{b}}(u)\in\mathfrak{b}$, $c_0\in \mathfrak{b}_{\bar 0}$ and $\psi(u),t,s\in\K$. Since $\omega_{\mathfrak{b}}$ is non-degenerate, there exists $a_{0}\in\mathfrak{b}_{\bar 0}$ such that
$
\psi(u)=\omega_{\mathfrak{b}}(a_{0},u),$ for all $u\in\mathfrak{b}.
$
The $\omega$-antisymmetry of $\rho$ implies $s=-\lambda$ and $c_{0}=a_{0}$. Moreover, since $\rho$ is $\omega$-antisymmetric, the map $\rho_{\mathfrak{b}}$ is $\omega_{\mathfrak{b}}$-antisymmetric as well, and because $\rho$ is invertible, it follows that  
$\rho_{\mathfrak{b}}$ is  invertible as well. Using the fact that $\rho\circ\mathrm{ad}_u$ is $\omega$-symmetric for all $u\in\mathfrak{b}$,  
we deduce that $\rho_{\mathfrak{b}}\circ\mathrm{ad}_{u}$ is $\omega_{\mathfrak{b}}$-symmetric. Therefore, the bilinear form
$
{\mathscr B}_{\mathfrak{b}}(-,-)
   =\omega_{\mathfrak{b}}\bigl(\rho_{\mathfrak{b}}(-),-\bigr)
$
is an even quadratic form on $\mathfrak{b}$,  
and 
$
(\mathfrak{b},\omega_{\mathfrak{b}}, \rho_{\mathfrak{b}})
$
is a flat QQF-superalgebra, such that $\omega_{\mathfrak b}$ is even.

Arguing as in the proof of Theorem~\ref{THQ1}, one checks that
\[
\rho_{\mathfrak{b}} \circ ( 2\xi-\xi^{*} )+\lambda(\xi+\xi^*) = \Rr_{a_0}^{\star_\mathfrak{b}}+ (\Rr_{a_0}^{\star_\mathfrak{b}})^*, 
\qquad \text{and} \qquad (2\xi^* - \xi)(a_0) = \lambda\, b_0, 
\]
 and therefore
$
(\g,\omega,\rho)
$
is an even flat double extension of 
$
(\mathfrak{b},\omega_{\mathfrak{b}},\rho_{\mathfrak{b}})
$
by means of $(\xi,b_{0}, a_0, \lambda)$.

If $e$ is odd, then
$
\omega(e,e)={\mathscr B}(\rho(e),e)=\lambda {\mathscr B}(e,e)=0,
$ since $\mathscr B$ is  symmetric. Hence $I:=\mathbb{K}e$ is totally isotropic.
The rest of the proof follows analogously.
\end{proof}

\sssbegin{Corollary}\label{suite-ex}
A  QQF-superalgebra $(\g, \omega,  \rho)$, for which  $|\rho|=|\omega|=\bar 0$, is flat if and only if it can be obtained by a sequence 
of even or odd flat double extensions of flat QQF-superalgebras, starting from the trivial superalgebra $\{0\}$.
\end{Corollary}
\begin{proof}
 Theorem \ref{quadratic-ex1} implies that $(\g, \om, \rho)$ is an even or odd flat double extension of a flat QQF-superalgebra
$(\mathfrak{b}_1,  \om_{\mathfrak{b}_1}, \rho_{\mathfrak{b}_1})$. The latter is itself an even or an odd flat double extension of a flat QQF-superalgebra
$(\mathfrak{b}_2,  \om_{\mathfrak{b}_2}, \rho_{\mathfrak{b}_2})$. Since the total dimension of $\g$ is finite and even (see Prop. \ref{dim-even}), there exists $k\in \mathbb{N}^*$ such that $(\mathfrak{b}_k,  \om_{\mathfrak{b}_k}, \rho_{\mathfrak{b}_k})$ is a trivial Lie superalgebra, i.e., $\mathfrak{b}_k=\{0\}$.
\end{proof}

We now establish the converses of Theorems \ref{THQ1} and \ref{THQ2} in the case where $|\omega|=\bar 1$. 
\sssbegin{Theorem}\label{quadratic-ex2}
Let $(\g, \om, \rho)$ be a flat QQF-superalgebra, for which $|\rho|=|\omega|+{\bar 1}=\bar 0$. Then $(\g, \om, \rho)$  is either:
\begin{enumerate}
    \item[$(i)$] an even flat  double extension of a flat QQF-superalgebra $(\mathfrak{b}, \om_\mathfrak{b}, \rho_\mathfrak{b})$ by means of $(\xi, b_0, a_0, \lambda)$, or
    \item[$(ii)$] an odd flat  double extension of a flat QQF-superalgebra $(\mathfrak{b}, \om_\mathfrak{b}, \rho_\mathfrak{b})$ by means of \allowbreak $(\xi, b_0, a_1, \lambda)$.
\end{enumerate}
\end{Theorem}

\begin{proof}
The proof is similar to that of Theorem \ref{quadratic-ex1}.
\end{proof}

\sssbegin{Corollary}\label{lcro}
A QQF-superalgebra $(\g, \om, \rho)$, for which $|\rho|=|\omega|+\bar 1=\bar 0$, is flat if and only if it can be obtained by a sequence 
of even or odd flat double extensions of flat  QQF-superalgebras, starting from the trivial superalgebra $\{0\}$.
\end{Corollary}

\begin{proof}
The proof is similar to that of Corollary \ref{suite-ex}.
\end{proof}
 \subsection{The case where \texorpdfstring{$|\rho|=\bar 1$}{|rho|=1}.}\label{oddrho}
We introduce here the concept of double extension for flat QQF–superalgebras in the case where $\rho$ is odd; that is, when the quadratic and quasi-Frobenius structures possess different parities. In this setting, the standard construction of a double extension by a one-dimensional space is no longer adequate. To overcome this obstacle, we develop the corresponding double extension by a two-dimensional space.

Let $(\mathfrak{b}, \omega_{\mathfrak{b}})$ be a flat orthosymplectic quasi-Frobenius Lie superalgebra,  
 and let $\star_{\mathfrak{b}}$ be the natural symplectic product associated with it.  Assume that there exist an even linear map 
$\xi_0: \mathfrak{b} \to \mathfrak{b}$, and an odd linear map $\xi_1 : \mathfrak{b} \to \mathfrak{b}$,  
 even elements $b_0, c_1 \in \mathfrak{b}_{\bar{0}}$, 
 odd elements $b_1, c_0 \in \mathfrak{b}_{\bar{1}}$, 
satisfying the system
\begin{equation}\label{plane1}
\begin{aligned}
(\xi_0-\xi_0^*)(b_0) &= 0,&
(2\xi_1+\xi_1^*)(c_0) &= 0,\\
(\xi_1+\xi_1^*)(b_0)&=0,&
3\xi_0^*(c_0)-\xi_0(c_0-b_1)+\xi_1^*(b_0) &= 0,\\
(\xi_1+2\xi_1^*)(b_0)-2\xi_0(2c_0+b_1) &= 0,&
\xi_1(c_0-b_1)+\xi_1^*(4c_0-b_1)-3\xi_0^*(c_1) &= 0,\\
\xi_0(2b_1+c_0)+\xi_0^*(2c_0+b_1)
 &=0,&
\xi_1(4c_0+5b_1)+\xi_1^*(2b_1+c_0)
 +(\xi_0^*-\xi_0)(c_1) &= 0,\\
3\xi_0(c_1)-(\xi_1+\xi_1^*)(2c_0+b_1) &= 0,
\end{aligned}
\end{equation}
\begin{equation}\label{plane2}
\begin{aligned}
 \relax[\xi_1,\xi_1^*]
  &= -3\xi_1^2 + \Rr_{c_1}^{\star_\mathfrak b},& 
\xi_0^*\circ\xi_0
  &= \tfrac13\bigl(\Rr_{b_0}^{\star_\mathfrak b}
      +(\Rr_{b_0}^{\star_\mathfrak b})^*\bigr),\\    
      (\xi_0^*\circ\xi_1-\xi_1^*\circ\xi_0)(u)
  &= \tfrac13(-1)^{|u|}
     \bigl(\Rr_{c_0-b_1}^{\star_\mathfrak b}
           +(\Rr_{c_0-b_1}^{\star_\mathfrak b})^*\bigr)(u),& [\xi_0,\xi_0^*]
  &= \xi_0^2 - \tfrac13 \Rr_{b_0}^{\star_\mathfrak b},\\
([\xi_1,\xi_0^*-\xi_0]-\xi_1\circ\xi_0)(u)
  &= \tfrac13(-1)^{|u|}\Rr_{2b_1+c_0}^{\star_\mathfrak b}(u),
\end{aligned}
\end{equation}
\begin{equation}\label{plane3}
\begin{aligned}
\Ll_{c_1}^{\star_\mathfrak b}
  &= -(\xi_1+\xi_1^*)^2,& 
\Ll_{c_0+b_1}^{\star_\mathfrak b}(u)
  &= (-1)^{|u|}[\xi_0^*-\xi_0,\xi_1+\xi_1^*](u),
\end{aligned}
\end{equation}
\begin{equation}\label{plane4}
\begin{aligned}
\xi_0([u,v]_\mathfrak b)
 &= \mathrm L_u^{\star_\mathfrak b}(\xi_0(v))
    -(-1)^{|u||v|}
     \mathrm L_v^{\star_\mathfrak b}(\xi_0(u)),\\
\xi_1([u,v]_\mathfrak b)
 &= \mathrm L_u^{\star_\mathfrak b}(\xi_1(v))
    -(-1)^{|u||v|}
     \mathrm L_v^{\star_\mathfrak b}(\xi_1(u)),\\
     D_0(u\star_\mathfrak b v)
 &= D_0(u)\star_\mathfrak b v
    + u\star_\mathfrak b D_0(v)
    - \xi_0(u)\star_\mathfrak b v,\\
D_1(u\star_\mathfrak b v)
 &= D_1(u)\star_\mathfrak b v
    + (-1)^{|u|}u\star_\mathfrak b D_1(v)
    + \xi_1(u)\star_\mathfrak b v,
\end{aligned}
\end{equation}
and 
 \begin{equation}\label{plane5}
\begin{array}{rcl} 3\om_\mathfrak{b}(c_1, b_0)-3\om_\mathfrak{b}(c_0-b_1, c_0+b_1) & = &7 \om_\mathfrak{b}(2b_1+c_0, 2c_0+b_1),\\[2mm] 6\om_\mathfrak{b}(c_1, b_0) & = & \om_\mathfrak{b}(2c_0+b_1, 2c_0+b_1),
\end{array}
 \end{equation}
where we have put $D_0(u)=(\xi^*_0-\xi_0)(u)$ and $D_1(u) =- (-1)^{|u|}(\xi_1^* + \xi_1)(u)$, for all $u\in \mathfrak{b}$.

 \sssbegin{Theorem} \label{dext-plane1}
 Let $(\mathfrak{b}, \omega_{\mathfrak{b}})$ be a flat orthosymplectic quasi-Frobenius Lie superalgebra,  
 and let $\star_{\mathfrak{b}}$ be the natural symplectic product associated with it. Let $V = \mathbb{K}d_0 \oplus \mathbb{K}d_1$ be a two-dimensional vector superspace, where 
 $|d_0|=\bar 0$ and $|d_1|=\bar 1$, and let $V^* = \mathbb{K}e_0 \oplus \mathbb{K}e_1$ be its dual. Assume that there exist an even linear map 
$\xi_0: \mathfrak{b} \to \mathfrak{b}$, and an odd linear map $\xi_1 : \mathfrak{b} \to \mathfrak{b}$,  
 even elements $b_0, c_1 \in \mathfrak{b}_{\bar{0}}$, 
 odd elements $b_1, c_0 \in \mathfrak{b}_{\bar{1}}$, 
 and a scalar $T \in \mathbb{K}$ 
satisfying the system (\ref{plane1})--(\ref{plane5}).

 Define the vector superspace $
 \mathfrak{g} := \mathbb{K}d_0 \oplus \mathbb{K}d_1 \oplus \mathfrak{b} \oplus \mathbb{K}e_0 \oplus \mathbb{K}e_1,
$
 equipped with the Lie bracket
\begin{equation*}\label{crochet-plan}
 \begin{array}{lcllcl}
 [d_0,u]  & = &  \omega_{\mathfrak{b}}(b_0, u)e_0 + \omega_{\mathfrak{b}}(b_1, u)e_1+(\xi_0^* - 2\xi_0)(u), & [d_0, d_1] & = & -(c_0+b_1) + T e_1, \\[0.3em]
 [d_1,u] & = & \omega_{\mathfrak{b}}(c_0, u)e_0 + \omega_{\mathfrak{b}}(c_1, u)e_1 - (-1)^{|u|}(\xi^*_1 + 2\xi_1)(u), &[d_1, d_1] & = &  2c_1 - 2T e_0. \\[2mm]
 [u,v] & = & [u,v]_{\mathfrak{b}} 
 + \omega_{\mathfrak{b}}\big((\xi_0 + \xi_0^*)(u), v\big)e_0 \\[2mm]
 &&
 + \Big((-1)^{|v|}\omega_{\mathfrak{b}}(\xi_1(u), v) 
 + (-1)^{|u|}\omega_{\mathfrak{b}}(\xi_1^*(u), v)\Big)e_1, \\[0.3em]
 \end{array}
 \end{equation*}
for all $u,v \in \mathfrak{b}$,  and the even bilinear form: \textup{(}for $i=\bar 0, \bar 1$\textup{)}
 \[
 \omega|_{\mathfrak{b} \times \mathfrak{b}} = \omega_{\mathfrak{b}}, 
 \quad \omega(d_i,\mathfrak{b}) = \omega(e_i,\mathfrak{b}) = 0,\quad \omega(e_i,d_i) = -(-1)^{|i|}\omega(d_i,e_i) = 1.
 \]
Then $(\mathfrak{g}, \omega)$ is a flat orthosymplectic quasi-Frobenius Lie superalgebra.

 Moreover, the natural symplectic product $\star$ associated with $(\mathfrak{g}, \omega)$ is given by: \textup{(}zero products are omitted\textup{)}
\begin{equation}\label{natural-plan}
 \begin{aligned}
 d_0 \star u & = \tfrac{1}{3}\omega_{\mathfrak{b}}(b_0, u)e_0+ \tfrac{1}{3}\omega_{\mathfrak{b}}(2b_1+c_0, u)e_1\\
 & +(\xi_0^* - \xi_0)(u) , &  d_0\star d_1  & =  -\frac{1}{3}(2b_1+c_0),\\
 d_1 \star u & = \tfrac{1}{3}\omega_{\mathfrak{b}}(2c_0+b_1, u)e_0 +  \omega_{\mathfrak{b}}(c_1, u)e_1\\
 &-(-1)^{|u|}(\xi_1^* + \xi_1)(u),&  d_0\star d_0  & = \frac{1}{3}b_0, \\
 u \star d_0 & = \xi_0(u) -\tfrac{2}{3}\omega_{\mathfrak{b}}(b_0, u)e_0+ \tfrac{1}{3}\omega_{\mathfrak{b}}(c_0-b_1, u)e_1,&  d_1\star d_0 & =  \frac{1}{3}(2c_0+b_1)-T e_1, \\
u \star d_1 & =  \xi_1(u) +\tfrac{1}{3}\omega_{\mathfrak{b}}(c_0-b_1, u)e_0, & d_1\star d_1 & = c_1-T e_0,\\
u\star v  & = u\star_\mathfrak{b} v+\om_\mathfrak{b}(\xi_0(u), v)e_0\\
&+(-1)^{|v|} \om_\mathfrak{b}(\xi_1(u), v)e_1,
 \end{aligned}
 \end{equation}
 for all $u,v \in \mathfrak{b}$. 

 \end{Theorem}

The flat orthosymplectic quasi-Frobenius Lie superalgebra $(\g, \omega)$ is called the flat planar double extension of $(\mathfrak{b}, \omega_{\mathfrak{b}})$ by means of $(\xi_0, \xi_1, b_0, b_1, c_0, c_1)$.

 \begin{proof}
 It is straightforward to verify that $\omega$ is closed on $\mathfrak{g}$.  
According to Eq.~\eqref{cyclic13}, the symplectic product associated with  
 $(\mathfrak{g}, \omega)$ is determined by the Eq.~\eqref{natural-plan}.  
 A direct computation then shows that this product indeed defines a 
 left-symmetric superalgebra.  
 Consequently, $(\mathfrak{g}, \omega)$ is a flat orthosymplectic 
 quasi-Frobenius Lie superalgebra.
 \end{proof}
 
We now turn to adapting this method to flat planar double extensions of flat orthosymplectic quasi-Frobenius Lie superalgebras that additionally carry a quadratic structure of odd parity.
\sssbegin{Theorem} \label{THRE1}
Let $(\mathfrak{b}, \omega_{\mathfrak{b}}, \rho_{\mathfrak{b}})$ be a flat QQF-superalgebra, where $|\rho_{\mathfrak{b}}|=|\omega_{\mathfrak{b}}|+{\bar 1}=\bar 1$. Let $(\mathfrak{g}, \omega)$ be the even flat planar  double extension of  
$(\mathfrak{b}, \omega_{\mathfrak{b}})$ by means of $(\xi_0,\xi_1,b_0,c_1, b_1, c_0)$ as in Theorem \ref{dext-plane1}.  
 Suppose that there exist $a_0 \in \mathfrak{b}_{\bar{0}}$, $a_1 \in \mathfrak{b}_{\bar{1}}$ and $\lambda \in \mathbb{K}\backslash \{0\}$ such that 
 \begin{equation}\label{sy-rho}
\begin{aligned}
\rho_{\mathfrak  b} \circ (\xi_0^*-2\xi_0)(u)
 + \lambda(-1)^{|u|}(\xi_1^*-\xi_1)(u)
 &= -\bigl(\Rr_{a_1}^{\star_\mathfrak b}
          +(\Rr_{a_1}^{\star_\mathfrak b})^*\bigr)(u),\\
(-1)^{|u|}\rho_\mathfrak b\circ(\xi_1^*+2\xi_1)(u)
 + \lambda(\xi_0^*+\xi_0)(u)
 &= \bigl(\Rr_{a_0}^{\star_\mathfrak b}
          +(\Rr_{a_0}^{\star_\mathfrak b})^*\bigr)(u),\\
(2\xi_1^*+\xi_1)(a_1)-\rho_\mathfrak b(c_0+b_1)
 &= -\lambda c_1,\\
(2\xi_0^*-\xi_0)(a_0)-\rho_\mathfrak b(c_0+b_1)
 &= \lambda b_0,\\
(2\xi_1^*+\xi_1)(a_0)-2\rho_\mathfrak b(c_1)
 &= -\lambda c_0,\\
(\xi_0-2\xi_0^*)(a_1)
 &= \lambda b_1,\\
(2\xi_0^*-\xi_0)(a_0)
 -(2\xi_1^*+\xi_1)(a_1)
 &= \lambda(c_1+b_0),\\
\omega_\mathfrak b(a_0,c_1)
 &= \omega_\mathfrak b(a_1,c_0+b_1)
 = \lambda T.
\end{aligned}
\end{equation}
For any $t \in \mathbb{K}$, define the odd linear endomorphism $\rho$ of $\mathfrak{g}$ by: \textup{(}for all $u \in \mathfrak{b}$\textup{)}
 \begin{equation}\label{map-rho}
 \begin{array}{lcllcl}
 \rho(e_0) & = &  \lambda e_1, &  \rho(e_1) & = & \lambda e_0,\, \\[2mm] 
  \rho(d_1) & = & 
a_0 -t e_0  + \lambda d_0, &   \rho(d_0) & = &  
a_1+ t e_1 - \lambda d_1,\\[2mm]
 \rho(u) & = &  \rho_{\mathfrak{b}}(u) + \omega_{\mathfrak{b}}(a_1, u) e_0+\omega_{\mathfrak{b}}(a_0, u)\, e_1.
 \end{array}
 \end{equation}
Then $(\mathfrak{g}, \omega, \rho)$ is a flat  QQF-superalgebra, where $|\rho|=|\omega|+{\bar 1}=\bar 1$.
 \end{Theorem}

The flat QQF-superalgebra $(\g, \omega)$ is called the flat  planar double extension of $(\mathfrak{b}, \omega_{\mathfrak{b}})$ by means of $(\xi_0, \xi_1, b_0, b_1, c_0, c_1, a_0, a_1)$.
 \begin{proof}
 The proof is similar to that of Theorem \ref{THQ1}.
 \end{proof}
 We now prove the converse of Theorem \ref{THRE1}.
\sssbegin{Theorem}\label{thm-rho-od1}
 Let $(\g, \om, \rho)$ be a flat QQF-superalgebra, where $|\rho|=|\omega|+{\bar 1}=\bar 1$. Then $(\g, \om, \rho)$ is a  flat  planar double extension of a flat QQF-superalgebra $(\mathfrak{b}, \om_\mathfrak{b}, \rho_\mathfrak{b})$ by means of $(\xi_0,\xi_1, b_0, b_1, c_0, c_1,a_0,a_1)$.
 \end{Theorem}
 \begin{proof}
 By \cite[Theorem 6.1]{BE}, we  have $Z(\g)\neq 0$. Hence, by 
Prop.~\ref{roh-inv}, we obtain
 $
 \rho\bigl(Z(\g)\bigr)
    = Z(\g).
 $
 Since $\K$ is algebraically closed, there exists 
 $e\in Z(\g)$ such that $\rho(e)=\lambda e$ with 
 $\lambda\neq 0$. Let us write $e=e_0+e_1$, where 
 $e_0\in Z(\g)_{\bar 0}$ and 
 $e_1\in Z(\g)_{\bar 1}$. 
 Since $\rho$ is odd, it follows that 
 $
 \rho(e_0)=\lambda e_1$, and   $
 \rho(e_1)=\lambda e_0.
 $ Set \(I:=\K e_0 \oplus \K e_1\).  
 Since \(e_0, e_1 \in Z(\g)\), we have
 $
 e_i \star u = u \star e_i = 0,$ for all $u\in \g,$ and $ i=0,1.$ Therefore,  \(I\) is a two-sided graded ideal of \((\g,\star)\). Moreover, since 
 $
 \omega(u\star v,\, I)=0,$  for all\, $u,v\in \g,
 $
 it follows that \(I^{\perp_\omega}\) is also a two-sided ideal of 
 \((\g,\star)\).  
 Because \(\omega\) is even, the ideal \(I^{\perp_\omega}\) is graded. Moreover, since $\omega$ is even and $\rho$ is $\om$-anti-symmetric, we have 
\[
 \omega(e_0,e_0)=0,\qquad 
 \lambda\,\omega(e_1,e_1)
  =\omega(\rho(e_0),e_1)
  =-\omega(e_0,\rho(e_1))
  =-\lambda\,\omega(e_0,e_0)=0,
 \]
 and
 \[
 \lambda\,\omega(e_0,e_1)
  =\omega(\rho(e_1),e_1)
  =\omega(e_1,\rho(e_1))
  =\lambda\,\omega(e_1,e_0)
  =- \lambda \omega(e_0,e_1)=0.
 \]
 Hence $I\subseteq I^\perp$.  Since $\omega$ is even and non-degenerate, there exist elements  
$d_0\in\g_{\bar{0}}\setminus\{0\}$ and  
 $d_1\in\g_{\bar{1}}\setminus\{0\}$ such that 
 \[
 \omega(e_0,d_0)=-\omega(d_0,e_0)=1,\qquad
 \omega(e_1,d_1)=\omega(d_1,e_1)=1,
 \]
 and
 \[
 \omega(e_0,d_1)=\omega(e_1,d_0)=\omega(d_0,d_1)=0.
 \]
 Since $I\subseteq I^\perp$, there exists a subsuperspace $\mathfrak{b}$ such that  $I^\perp = I \oplus \mathfrak{b}$ and  the restriction  
 $\omega_{\mathfrak{b}}=\omega|_{\mathfrak{b}\times \mathfrak{b}}$  
 is non-degenerate.  We may therefore write
 \(
 \mathfrak{g}
 = \mathbb{K}e_0 \oplus \mathbb{K}e_1 \oplus
 \mathfrak{b} \oplus 
 \mathbb{K}d_0 \oplus \mathbb{K}d_1.
 \)

Since $I^{\perp}=I\oplus\mathfrak{b}$ is a two-sided graded ideal of 
 $(\mathfrak{g},\star)$, for any  $u,v\in\mathfrak{b}$ we have
\[
u\star v
 = u\star_{\mathfrak{b}} v 
   + \mu(u,v)\,e_0
   + \psi(u,v)\,e_1,
 \]
 where 
 $\mu,\psi:\mathfrak{b}\times\mathfrak{b}\to\mathbb{K}$ 
 are even and odd bilinear maps respectively, and 
$\star_{\mathfrak{b}}:\mathfrak{b}\times\mathfrak{b}\to\mathfrak{b}$ 
 is an even bilinear product.

 It follows that $({\mathfrak b},\star_{\mathfrak b})$ is a left-symmetric
 superalgebra, that $\mu$ and $\psi$ satisfy the compatibility conditions
$$
\begin{array}{lcl}\mu([u,v], w)& = & \mu(u\star_\mathfrak{b} v, w)-(-1)^{|u||v|}\mu(v\star_\mathfrak{b} u, w),\\[2mm] 
\psi([u,v], w)&=&\psi(u\star_\mathfrak{b} v, w)-(-1)^{|u||v|}\psi(v\star_\mathfrak{b} u, w),
\end{array}
$$
for all $u,v,w\in \mathfrak{b}$, and that 
$\omega_{\mathfrak{b}}$ is orthosymplectic and closed on
 ${\mathfrak{b}}$. Since $I^{\perp}$ is also a graded ideal of $\mathfrak{g}$,
 the Lie brackets on $\mathfrak{g}$ take the form 
\[
 \begin{array}{lcl}
 [d_0,u] &=& D_0(u) + L_0(u)e_0 + L_1(u)e_1,\\[2mm]
[d_1,u] &= &D_1(u) + S_0(u)e_0 + S_1(u)e_1,\\[2mm]
[u,v] &= &[u,v]_{\mathfrak{b}}
 + \big(\mu(u,v)-(-1)^{|u||v|}\mu(v,u)\big)e_0
 + \big(\psi(u,v)-(-1)^{|u||v|}\psi(v,u)\big)e_1,\\[2mm]
 [d_0,d_1] &=& T e_1 + x_1 + \beta_1 d_1,\\[2mm]
[d_1,d_1] &= &T_0 e_0 + x_0 + \beta_0 d_0,
\end{array}
\]
 where 
 $
 D_0 \in \mathrm{End}(\mathfrak{b})_{\bar{0}}, 
 D_1 \in \mathrm{End}(\mathfrak{b})_{\bar{1}},
 $ the maps $L_0,S_1:\mathfrak{b}\to\mathbb{K}$ are even, the maps  $L_1,S_0:\mathfrak{b}\to\mathbb{K}$ are odd, 
 $x_0\in\mathfrak{b}_{\bar{0}}$, $ x_1\in\mathfrak{b}_{\bar{1}},$ and $
 \alpha_0,\alpha_1,\beta_0,\beta_1\in\mathbb{K}.$

 Since $\omega_{\mathfrak{b}}$ is non-degenerate, there exist 
$b_0,c_1\in\mathfrak{b}_{\bar{0}}$ and $c_0,b_1\in\mathfrak{b}_{\bar{1}}$ such that
 \[
 L_0=\omega_{\mathfrak{b}}(b_0,\cdot),\qquad
 L_1=\omega_{\mathfrak{b}}(b_1,\cdot),\qquad
 S_0=\omega_{\mathfrak{b}}(c_0,\cdot),\qquad
 S_1=\omega_{\mathfrak{b}}(c_1,\cdot).
\]
Since $\omega$ is closed on $\mathfrak{g}$, it follows that 
 \[
 x_1=-b_1-c_0,\qquad 
 x_0=2c_1,\qquad 
 T_0=-2T,\qquad
 \beta_0=\beta_1=0.
 \]
Moreover, the product associated with \( (\g, \omega) \) is given by
\[
\begin{array}{lcllcl}
u \star v &=& u \star_\mathfrak{b} v + \mu(u,v)e_0+\psi(u,v) e_1,& 
d_0 \star u & = & G_0(u) + h_0(u)e_0+h_1(u)e_1, \\[2mm]
d_1 \star u & = &  G_1(u) + k_0(u)e_0+k_1(u)e_1, & 
u \star d_0 & = & \xi_0(u) + f_0(u)e_0+f_1(u) e_1, \\[2mm]
u \star d_1 & = & \xi_1(u) + g_0(u)e_0+g_1(u) e_1, & 
 d_0 \star d_0 & = & \lambda_0 d_0 + m_0 + t_0 e_0,\\[2mm]
 d_0 \star d_1 & = &  \lambda_1 d_1 + m_1 + t_1 e_1,& 
 d_1 \star d_0 & = &  \lambda_2 d_1 + m_1' + t_2 e_1,\\[2mm] 
 d_1 \star d_1 & = &  \lambda_3 d_0 + m_0' + t_3 e_0,
\end{array}
\]
where \( \xi_0, G_0 \in \mathrm{End}(\mathfrak{b})_{\bar{0}} \), \( \xi_1, G_1 \in \mathrm{End}(\mathfrak{b})_{\bar{1}} \),   \( f_0, g_1, h_0, k_1 : \mathfrak{b} \to \K \) are even linear maps, \( f_1, g_0, h_1, k_0 : \mathfrak{b} \to \K \) are odd linear maps, and $m_0, m_0'\in \mathfrak{b}_{\bar{0}}$, $m_1, m_1'\in \mathfrak{b}_{\bar{1}}$ and $\la_i, t_i\in \K$, with $i\in\{0,1,2,3\}$.

By the natural symplectic product \eqref{cyclic13}, we have for any $u,v\in \mathfrak{b}$: 
\[
\mu(u,v) = \omega(u \star v, d_0) = -(-1)^{|u||v|}\omega(v, u \star d_0)
= \omega_\mathfrak{b}(\xi_0(u), v),
\]
and hence
\[
\omega_\mathfrak{b}(\xi_0(u), v)=\om(u\star v, d_0) =\tfrac{1}{3}\omega ([u, v],d_0)+\tfrac{1}{3}\om([u,d_0], v)= \tfrac{1}{3}\omega_\mathfrak{b}((\xi_0 + \xi_0^* - D_0)(u), v),
\]
which gives \( D_0 = \xi_0^* - 2\xi_0 \).
Similarly,  for any $u,v\in \mathfrak{b}$, we have
\[
\psi(u,v) = \omega(u \star v, d_1) =- (-1)^{|u||v|}\omega(v, u \star d_1)
= (-1)^{|v|}\omega_\mathfrak{b}(\xi_1(u), v),
\]
and hence 
\begin{align*}
(-1)^{|v|}\,\omega_{\mathfrak{b}}(\xi_1(u), v)
&= \omega(u \star v, d_1)
= \tfrac{1}{3}\om([u, v], d_1) + \tfrac{1}{3}(-1)^{|v|}\omega([u, d_1], v) \\
&= \tfrac{1}{3}(-1)^{|v|}\omega_{\mathfrak{b}}(\xi_1(u), v)
+ \tfrac{1}{3}(-1)^{|u|}\omega_{\mathfrak{b}}(\xi_1^{*}(u), v)
- \tfrac{1}{3}(-1)^{|u|+|v|}\omega_{\mathfrak{b}}(D_1(u), v).
\end{align*}
Hence,
\[
\omega_{\mathfrak{b}}(D_1(u), v)
= -2(-1)^{|u|}\omega_{\mathfrak{b}}(\xi_1(u), v)
+ (-1)^{|u|+|v|}\omega_{\mathfrak{b}}(\xi_1^{*}(u), v).
\]
Since both $\xi_1$ and $D_1$ are odd, and $\omega$ is even, it follows that when $|u| + |v| = \bar{0}$, the equality reduces to the trivial identity $0 = 0$. Consequently, the nontrivial case occurs when $|u| + |v| = \bar{1}$, leading to the following relation
\[
\omega_{\mathfrak{b}}(D_1(u), v)
=- (-1)^{|u|}\omega_{\mathfrak{b}}((2\xi_1+\xi_1^*)(u), v).
\]
It follows that $D_1(u)=- (-1)^{|u|}(2\xi_1+\xi_1^*)(u)$. In the same way, one shows that
\[
\begin{array}{lcllcllcllcl}
G_0(u) & = & (\xi_0^{*} - \xi_0)(u), &
G_1(u) & = & (-1)^{|u|}\,(\xi_1^{*} + \xi_1)(u), &
h_0(u) & = & \tfrac{1}{3}\,\omega_{\mathfrak{b}}(b_0, u), \\[2mm] 
h_1(u) & = & \tfrac{1}{3}\,\omega_{\mathfrak{b}}(2b_1 + c_0, u), &k_0(u) & = & \tfrac{1}{3}\,\omega_{\mathfrak{b}}(2c_0 + b_1, u), & 
k_1(u) & = & \omega_{\mathfrak{b}}(c_1, u), \\[2mm]
f_0(u) & = & \tfrac{2}{3}\,\omega_{\mathfrak{b}}(b_0, u), & 
f_1(u) & = & \tfrac{1}{3}\,\omega_{\mathfrak{b}}(c_0 - b_1, u),& g_0(u) & = &  \tfrac{1}{3}\,\omega_{\mathfrak{b}}(c_0 - b_1, u), \\[2mm]
g_1(u) & = & 0, & 
m_0 & = & \tfrac{1}{3} b_0, & 
m_1 & = & -\tfrac{1}{3}(2b_1 + c_0), \\[2mm] 
m_1' & = &  \tfrac{1}{3}(2c_0 + b_1), & 
m_0' & = & c_1,
\end{array}
\]
and 
\[
\lambda_0 = \lambda_1 = \lambda_2 = \lambda_3 = 
t_0 = t_1 = 0, \qquad
t_2 = t_3 = -T.
\]

Consequently, the natural symplectic product associated with 
$(\g, \om)$ takes the form as in   \eqref{natural-plan}

 A direct computation shows that the natural symplectic product
 \eqref{natural-plan} is left-symmetric if and only if 
 $(\xi_0,\xi_1,b_0,, b_1, c_0, c_1)$ satisfies the system of equations \eqref{plane1}, \eqref{plane2}, \eqref{plane3}, \eqref{plane4}, \eqref{plane5}.
 Therefore, $
 (\mathfrak{g},\omega)$ is a flat planar double extension of $
 (\mathfrak{b},\omega_{\mathfrak{b}})$
 by means of $
 (\xi_0,\xi_1,b_0,, b_1, c_0, c_1).$

 Finally, since $\rho(I)=I$ and $\rho$ is $\omega$-antisymmetric, we have $\rho(I^\perp)=I^\perp$. Writing
 \[
 \begin{aligned}
 \rho(u) &= \rho_{\mathfrak{b}}(u) + f(u)e_0 + g(u)e_1, \quad \text{for all } u\in \mathfrak{b},\\
 \rho(d_0) &= t_1 e_1 + a_1 + s_1 d_1,\\
 \rho(d_1) &= t_0 e_0 + a_0 + s_0 d_0,
 \end{aligned}
 \]
 with $a_0\in\mathfrak{b}_{\bar 0}, a_1\in\mathfrak{b}_{\bar 1}$ and scalars $f(u),g(u),t_0,t_1,s_0,s_1\in\mathbb{K}$. The $\omega$-antisymmetry of $\rho$ gives
\[
 f(u) = \omega_{\mathfrak{b}}(a_1,u), \quad g(u) = \omega_{\mathfrak{b}}(a_0,u), \quad t_0=-t_1, \quad s_0=-s_1=\lambda.
\]
Hence, $\rho$ has the form given in \eqref{map-rho}. Moreover, $\rho_{\mathfrak{b}}$ is $\omega_{\mathfrak{b}}$-antisymmetric and invertible. Since $\rho\circ \mathrm{ad}_u$ is $\omega$-antisymmetric for all $u\in \mathfrak{b}$, it follows that $\rho_{\mathfrak{b}}\circ \mathrm{ad}_u$ is symmetric with respect to $\omega_{\mathfrak{b}}$. Therefore, the bilinear form
 $
 {\mathscr B}_{\mathfrak{b}}(-,-) = \omega_{\mathfrak{b}}(\rho_{\mathfrak{b}}(-),-)
$, defines an odd quadratic form on $(\mathfrak{b},\omega_{\mathfrak{b}})$, making 
$
(\mathfrak{b},\omega_{\mathfrak{b}},\rho_{\mathfrak{b}})
 $,  
a flat QQF-superalgebra. 

It is straightforward to check that $\rho_{\mathfrak{b}}$ satisfies the system \eqref{sy-rho}, so that 
 $
 (\g,\omega,\rho)
 $
 is a  flat  quadratic planar double extension of 
 $
 (\mathfrak{b},\omega_{\mathfrak{b}},\rho_{\mathfrak{b}})
 $
 by means of $(\xi_0,\xi_1,b_0, b_1, c_0, c_1)$,
 \end{proof}
 
 \sssbegin{Corollary}\label{cosuite}
Let $(\g, \om, \rho)$ be a QQF-superalgebra with $|\rho|=|\omega|+\bar 1= \bar 1$.  Then the total dimension of $\g$ is $4n$. Moreover,   $(\g, \om, \rho)$ is flat if and only if it can be obtained by a sequence 
 of flat planar double extensions starting from the trivial superalgebra $\{0\}$.

 \end{Corollary}

\begin{proof}
Theorem \ref{thm-rho-od1} implies that $(\g, \om, \rho)$ is a flat planar double extension of a QQF-superalgebra
$(\mathfrak{b}_1,  \om_{\mathfrak{b}_1}, \rho_{\mathfrak{b}_1})$ by means of  $(\xi_0, \xi_1, b_0, b_1, c_0,  c_1)$. The latter is itself a flat planar double extension of flat QQF-superalgebra
$(\mathfrak{b}_2,  \om_{\mathfrak{b}_2}, \rho_{\mathfrak{b}_2})$ by means of  $(\xi'_0, \xi'_1, b'_0, b'_1, c'_0, c'_1)$. Recall that the total dimension of $\g$ is even (see Prop. \ref{dim-even}). Suppose $\dim\g= 4 n$, where $n\in \mathbb{N}$. It follows that $\mathrm{dim}(\mathfrak{b}_1)=4(n-1).$ Repeating this procedure $n$  times, it follows that $(\mathfrak{b}_n,  \om_{\mathfrak{b}_n}, \rho_{\mathfrak{b}_n})$ is the trivial Lie superalgebra, i.e., $\mathfrak{b}_n=\{0\}$.

Likewise, if $\dim\g= 4 n+2$, where $n\in \mathbb{N}$, then $\mathrm{dim}({\mathfrak b}_1)=4(n-1)+2$. Repeating this procedure $n$  times, it follows that $(\mathfrak{b}_n,  \om_{\mathfrak{b}_n}, \rho_{\mathfrak{b}_n})$ is a 2-dimensional flat QQF-superalgebra. However, a 
$2$-dimensional Lie superalgebra cannot admit both a symplectic and a quadratic structure of different parities; this yields a contradiction. Therefore, the dimension of $\mathfrak{g}$ must always be $4n$. This completes the proof.
\end{proof}

Let $(\mathfrak{b}, \omega_{\mathfrak{b}})$ be a flat periplectic quasi-Frobenius Lie superalgebra,  
 and let $\star_{\mathfrak{b}}$ be the natural symplectic product associated with it.  Assume that there exist an even linear map 
 $\xi_0: \mathfrak{b} \to \mathfrak{b}$, an odd linear map $\xi_1 : \mathfrak{b} \to \mathfrak{b}$, 
 even elements $b_1, c_0 \in \mathfrak{b}_{\bar{0}}$, 
 odd elements $b_0, c_1 \in \mathfrak{b}_{\bar{1}}$, 
 and a scalar $T \in \mathbb{K}$ 
 satisfying the system 
 \begin{equation}\label{2plane1}
 \begin{aligned}
(\xi_0^*-\xi_0)(b_1) &= 0,&  
(\xi_1^*-2\xi_1)(c_0) &= 0,\\
\xi_0(b_0+c_1)-3\xi_0^*(c_1)+\xi_1^*(b_1) &= 0,&  
\xi_1(b_0+c_1)-\xi_1^*(4c_1+b_0)+3\xi_0^*(c_0) &= 0,\\
(\xi_1-2\xi_1^*)(b_1)+\xi_0^*(2c_1-b_1) &= 0,& 
\xi_0(2b_0-c_1)-\xi_0^*(2c_1-b_0)
 +(\xi_1^*-\xi_1)(b_1) &= 0,\\
\xi_1(5b_0-4c_1)-\xi_1^*(2b_0-c_1)
\\ +3(\xi_0^*-\xi_0)(c_0) &= 0,&
(\xi_1-\xi_1^*)(2c_1-b_0)+3\xi_0(c_0) &= 0,
\end{aligned}
\end{equation}
\begin{equation}\label{2plane2}
\begin{aligned}
\relax [\xi_0,\xi_0^*]
  &= \xi_0^2-\tfrac13 \Rr_{b_1}^{\star_\mathfrak b}, & (\xi_0^*\circ\xi_1+\xi_1^*\circ\xi_0)(u)
 &= \tfrac13(-1)^{|u|}
    \bigl(\Rr_{b_0+c_1}^{\star_\mathfrak b}
          +(\Rr_{b_0+c_1}^{\star_\mathfrak b})^*\bigr)(u),\\
[\xi_1,\xi_1^*]
  &= 3\xi_1^2+\Rr_{c_0}^{\star_\mathfrak b}, &[\xi_1,\xi_0^*-\xi_0](u)=&(\xi_1\circ\xi_0)(u)
 -\tfrac13(-1)^{|u|}
    \Rr_{2b_0-c_1}^{\star_\mathfrak b}(u),\\
\xi_0^*\circ\xi_0
  &= \tfrac13\bigl(\Rr_{b_1}^{\star_\mathfrak b}
     +(\Rr_{b_1}^{\star_\mathfrak b})^*\bigr), & [\xi_0,\xi_1^*-\xi_1](u)=&(\xi_0\circ\xi_1)(u) -\tfrac13(-1)^{|u|}
    \Rr_{2c_1-b_0}^{\star_\mathfrak b}(u),
\end{aligned}
\end{equation}
\begin{equation}\label{2plane3}
\begin{aligned}
\Ll_{c_0}^{\star_\mathfrak b}
 &= (\xi_1^*-\xi_1)^2,& 
\Ll_{b_0-c_1}^{\star_\mathfrak b}(u)
 &= (-1)^{|u|}
    [\xi_0^*-\xi_0,\xi_1^*-\xi_1](u),
    \end{aligned}
\end{equation}
\begin{equation}\label{2plane4}
\begin{aligned}
\xi_0([u,v]_\mathfrak b)
 &= \mathrm L_u^{\star_\mathfrak b}(\xi_0(v))
    -(-1)^{|u||v|}
     \mathrm L_v^{\star_\mathfrak b}(\xi_0(u)),\\
\xi_1([u,v]_\mathfrak b)
 &= \mathrm L_u^{\star_\mathfrak b}(\xi_1(v))
    -(-1)^{|u||v|}
     \mathrm L_v^{\star_\mathfrak b}(\xi_1(u)),\\
(\xi_0^*-\xi_0)(u\star_\mathfrak b v)
 &= (\xi_0^*-\xi_0)(u)\star_\mathfrak b v
    + u\star_\mathfrak b(\xi_0^*-\xi_0)(v)
    - \xi_0(u)\star_\mathfrak b v,\\
(-1)^{|u|+|v|}(\xi_1^*-\xi_1)(u\star_\mathfrak b v)
 &= (-1)^{|u|}(\xi_1^*-\xi_1)(u)\star_\mathfrak b v\\
 &\quad
    +(-1)^{|u|+|v|}
     u\star_\mathfrak b(\xi_1^*-\xi_1)(v)
    -(-1)^{|u|}\xi_1(u)\star_\mathfrak b v,
\end{aligned}
\end{equation}
\begin{equation}\label{2plane6}
\begin{aligned}
\omega_\mathfrak b(b_1,2b_0-c_1) &= 0, & 
\omega_\mathfrak b(c_0,c_1) &= 0,
\end{aligned}
\end{equation}
\sssbegin{Theorem}\label{THDOUBLE} Let $(\mathfrak{b}, \omega_{\mathfrak{b}})$ be a flat periplectic quasi-Frobenius Lie superalgebra,  
 and let $\star_{\mathfrak{b}}$ be the natural symplectic product associated with it. Let $V = \mathbb{K}d_0 \oplus \mathbb{K}d_1$ be a two-dimensional vector superspace, where 
 $V_{\bar{0}} = \mathbb{K}d_0$ and $V_{\bar{1}} = \mathbb{K}d_1$, 
 and let $\Pi(V^*) = \mathbb{K}e_0 \oplus \mathbb{K}e_1$, where $V^*$ is its dual. Assume that there exist an even linear map 
 $\xi_0: \mathfrak{b} \to \mathfrak{b}$, an odd linear map $\xi_1 : \mathfrak{b} \to \mathfrak{b}$, 
 even elements $b_1, c_0 \in \mathfrak{b}_{\bar{0}}$, 
 odd elements $b_0, c_1 \in \mathfrak{b}_{\bar{1}}$, 
 and a scalar $T \in \mathbb{K}$ 
 satisfying the system \eqref{2plane1}-- \eqref{2plane6}. 
Define the vector superspace $
\mathfrak{g} := \mathbb{K}d_0 \oplus \mathbb{K}d_1 \oplus \mathfrak{b} \oplus \mathbb{K}e_0 \oplus \mathbb{K}e_1,
$
equipped with the Lie bracket
 \begin{equation*}
 \begin{aligned}
\relax  [d_0,u] &= (\xi_0^* - 2\xi_0)(u) + \omega_{\mathfrak{b}}(b_0, u)e_0 + \omega_{\mathfrak{b}}(b_1, u)e_1, \\[0.3em]
 [d_1,u] &=  (-1)^{|u|}(\xi^*_1 - 2\xi_1)(u)+ \omega_{\mathfrak{b}}(c_0, u)e_0 + \omega_{\mathfrak{b}}(c_1, u)e_1, \\[0.3em]
 [u,v] &= [u,v]_{\mathfrak{b}} 
+  \Big((-1)^{|v|}\omega_{\mathfrak{b}}(\xi_1(u), v) 
 + (-1)^{|u|}\omega_{\mathfrak{b}}(\xi_1^*(u), v)\Big)e_0+\omega_{\mathfrak{b}}\big((\xi_0 + \xi_0^*)(u), v\big)e_1  , \\[0.3em]
 [d_0, d_1] & =b_0-c_1 + T e_1, \qquad [d_1, d_1] = -2c_0.
 \end{aligned}
 \end{equation*}
for all $u,v \in \mathfrak{b}$, and an even bilinear form $\omega$ given by
 \[
 \omega|_{\mathfrak{b} \times \mathfrak{b}} = \omega_{\mathfrak{b}}, 
 \quad \omega(d_i,\mathfrak{b}) = \omega(e_i,\mathfrak{b}) = 0,\quad \omega(e_i,d_j) = -\omega(d_j,e_i) = 1,
 \]
 for all $i, j \in \{0,1\}$, with $i\neq j$.

 Then $(\mathfrak{g}, \omega)$ is a flat periplectic  quasi-Frobenius Lie superalgebra.

Moreover, the natural symplectic product $\star$ associated with $(\mathfrak{g}, \omega)$ is given by: \textup{(}zero products are omitted\textup{)} 
 \begin{equation*}\label{natural-plan3}
 \begin{aligned}
 d_0 \star u & = (\xi_0^* - \xi_0)(u) + \tfrac{1}{3}\omega_{\mathfrak{b}}(2b_0-c_1, u)e_0+ \tfrac{1}{3}\omega_{\mathfrak{b}}(b_1, u)e_1,& d_1\star d_1 & =-c_0 \\[1mm]
 d_1 \star u & = (-1)^{|u|}(\xi_1^* - \xi_1)(u) +\omega_{\mathfrak{b}}(c_0, u)e_0+ \frac{1}{3}\omega_{\mathfrak{b}}(2c_1-b_0, u)e_1,& d_1\star d_0 & =\frac{1}{3}(2c_1-b_0)-\frac{2}{3}T e_1, \\[1mm]
 u \star d_0 & = \xi_0(u) -\tfrac{1}{3}\omega_{\mathfrak{b}}(b_0+c_1, u)e_0 -\tfrac{2}{3}\omega_{\mathfrak{b}}(b_1, u)e_1 &  d_0\star d_1 &=\frac{1}{3}(2b_0-c_1)+\frac{1}{3}T e_1,\\[1mm]
 u \star d_1 & = \xi_1(u)-\tfrac{1}{3}\omega_{\mathfrak{b}}(b_0+c_1, u)e_1, &  d_0\star d_0 &=\frac{1}{3}b_1,\\[1mm]
 u\star v & =u\star_\mathfrak{b} v+(-1)^{|v|} \om_\mathfrak{b}(\xi_1(u), v)e_0 +\om_\mathfrak{b}(\xi_0(u), v)e_1,
 \end{aligned}
 \end{equation*}
 for all $u,v \in \mathfrak{b}$.

 \end{Theorem}

The flat periplectic quasi-Frobenius superalgebra $(\g, \omega)$ is called the flat  planar double extension of $(\mathfrak{b}, \omega_{\mathfrak{b}})$ by means of $(\xi_0, \xi_1, b_0, b_1, c_0, c_1)$.

We now adapt this method to the case of flat planar double extensions of flat periplectic quasi-Frobenius Lie superalgebras endowed with a quadratic structure of odd parity.
\sssbegin{Theorem} \label{THRE2}
Let $(\mathfrak{b},  \omega_{\mathfrak{b}}, \rho_{\mathfrak{b}})$ be a flat QQF-superalgebra, where $|\rho_{\mathfrak{b}}|=|\omega_{\mathfrak{b}}|=\bar 1$. Let  
 $(\mathfrak{g}, \omega)$ be the  flat planar  double extension of  
$(\mathfrak{b}, \omega_{\mathfrak{b}})$ by means of $(\xi_0, \xi_1, b_0, b_1, c_0, c_1)$.  
 Suppose that there exist $a_0 \in \mathfrak{b}_{\bar{0}}$, $a_1 \in \mathfrak{b}_{\bar{1}}$ and $\lambda \in \mathbb{K}\backslash \{0\}$ such that 
\begin{align*}
(\xi_0-2\xi_0^*)(a_0)+\rho_\mathfrak{b}(b_0-c_1)
  &= \lambda b_1,
&\quad
(\xi_1-2\xi_1^*)(a_1)+\rho_\mathfrak{b}(c_1-b_0)
  &= \lambda c_0,\\[2mm]
(\xi_1-2\xi_1^*)(a_0)+2\rho_\mathfrak{b}(c_0)
  &= \lambda c_1,
&\quad
(2\xi_0^*-\xi_0)(a_1)
  &= \lambda b_0.
\end{align*}

\begin{align}\label{sy-rho1}
  (\xi_0-2\xi_0^*)(a_0)+(\xi_1-2\xi_1^*)(a_1) & =  \la(b_1-c_0),\\[2mm] 
(-1)^{|u|}\rho_\mathfrak{b}\circ (2\xi_0^*-\xi_0)(u)+\la(\xi_1+\xi^*_1)(u)& =  (-1)^{|u|}(\Rr^{\star_\mathfrak{b}}_{a_1}+(\Rr^{\star_\mathfrak{b}}_{a_1})^*)(u), \nonumber\\[2mm] \rho_\mathfrak{b}\circ (\xi_1-2\xi_1^*)(u)+(-1)^{|u|}\la(\xi_0+\xi^*_0)(u) & =  -(-1)^{|u|}(\Rr^{\star_\mathfrak{b}}_{a_0}+(\Rr^{\star_\mathfrak{b}}_{a_0})^*)(u),\nonumber\\
2\om_\mathfrak{b}(a_1, c_0)+\om_\mathfrak{b}(a_0, b_0-c_1)& =\la T.\nonumber
\end{align}
For any $t \in \mathbb{K}$, define the odd linear endomorphism $\rho$ of $\mathfrak{g}$ by: (for all $u \in \mathfrak{b}$)
 \begin{equation*}\label{map-rho2}
 \begin{cases}
 \rho(e_0) = \lambda e_1,\quad  \rho(e_1) = \lambda e_0,\quad  \rho(d_0) = 
 t e_1 + a_1 - \lambda d_1,\quad \rho(d_1) = 
  a_0 + \lambda d_0, \\ 
 \rho(u) = \rho_{\mathfrak{b}}(u) - \omega_{\mathfrak{b}}(a_0, u)\, e_0+\omega_{\mathfrak{b}}(a_1, u)\, e_1.
 \end{cases}
 \end{equation*}
Then $(\mathfrak{g}, \omega, \rho)$ is a flat QQF-superalgebra.
 \end{Theorem}

The flat QQF-superalgebra $(\g, \omega)$ is called the flat  planar double extension of $(\mathfrak{b}, \omega_{\mathfrak{b}})$ by means of $(\xi_0, \xi_1, b_0, b_1, c_0, c_1, a_0, a_1)$.

\sssbegin{Theorem}\label{thm-rho-od2}
 Let $(\g, \om, \rho)$ be a flat QQF-superalgebra, where $|\rho|=|\omega|=\bar 1$. Then $(\g, \om, \rho)$ is a  flat  planar double extension of a flat QQF-superalgebra $(\mathfrak{b}, \om_\mathfrak{b}, \rho_\mathfrak{b})$ by means of $(\xi_0,\xi_1, b_0, b_1, c_0, c_1, a_0, a_1)$.
\end{Theorem}

 \begin{proof}
The proof is similar to that of Theorem \ref{THRE1}.
\end{proof}

\sssbegin{Corollary}\label{cosuite1}
Let $(\g, \om, \rho)$ be a QQF-superalgebra with $|\rho|=|\omega|=\bar 1$.  Then, the total dimension of $\g$ is $4n$. Moreover,  $(\g, \om, \rho)$  is flat if and only if it can be obtained by a sequence 
 of  flat planar double extensions   starting from the trivial superalgebra $\{0\}$.
 \end{Corollary}

 \begin{proof}
The proof is similar to that of Corollary \ref{cosuite}.
\end{proof}
\section{Examples in low-dimension} \label{examples-low}
In this section, we classify flat non-abelian QQF-superalgebras of total dimension four and present examples in dimension six and eight. Recall that Proposition \ref{dim-even} implies that the total dimension is even, and by \cite{BE} any 2-dimensional flat quasi-Frobenius Lie superalgebra is necessarily abelian. 




 \ssbegin{Proposition}\label{rho-odd}
Any four-dimensional flat QQF-superalgebra $(\mathfrak{g}, \omega, \rho)$ with an odd endomorphism $\rho$ is necessarily abelian.
\end{Proposition}


\begin{proof}
Let $(\mathfrak{g}, \om, \rho)$ be  a flat non-abelian QQF-superalgebra of total dimension $4$ where $\rho$ is odd. Theorems~\ref{thm-rho-od1} and~\ref{thm-rho-od2}, imply that such a superalgebra is a  flat quadratic planar double extension of the  trivial  $\mathfrak{b}=\{0\}$. According to Theorems~\ref{THRE1} and~\ref{THRE2}, we get  $\xi_i=0$, $b_i=c_i=0$, for all $i\in\{0,1\}$, and $T=0$. Thus, $\g$ is an abelian Lie superalgebra.
\end{proof}

\ssbegin{Proposition}[\cite{BE}]\label{cls-dim4}
Any four-dimensional flat quasi-Frobenius non-abelian Lie superalgebra is isomorphic to one of the following Lie superalgebras:

\begin{itemize}
\item[$(i)$]  $(\mathbb{K}\oplus \h_3, \om)$ \textup{(}in the basis $\{x_1, x_2, x_3, x_4\}$\textup{)}:   
$$[x_1, x_2]=x_3, \quad \om= x_1^*\wedge x_4^*+x_2^* \wedge x_3^*,$$

\item[$(ii)$]  $(\g^2, \om)$ \textup{(}in the basis $\{x_1, x_2\mid  y_1, y_2\}$\textup{)}:
$$
	[x_1,y_1]=y_2,\quad [y_1, y_1]= x_2\mbox{ and }\om=2x_1^*\wedge x_2^*-y_1^*\wedge y_2^*,
	$$
\item[$(iii)$]  $(\g^3, \om)$ \textup{(}in the basis $\{x_1, x_2\mid  y_1, y_2\}$\textup{)}:
$$[x_1, y_1]=y_2, \quad \om= x_1^*\wedge y_2^*+x_2^* \wedge y_1^*,$$
\item[$(iv)$]  $(\g^4, \om)$ \textup{(}in the basis $\{x_1, x_2\mid  y_1, y_2\}$\textup{)}:
$$[y_1, y_1]=x_1,\quad [y_1, y_2]=x_2, \quad \om= -2x_1^*\wedge y_2^*+x_2^* \wedge y_1^*.$$
\end{itemize}
\end{Proposition}

\ssbegin{Proposition}
Any four-dimensional flat non-abelian QQF-superalgebra is isomorphic to one of the following Lie superalgebras:

\begin{itemize}
\item[$(i)$]  $(\g^2, \om, \rho)$ \textup{(}in the basis $\{x_1, x_2\mid  y_1, y_2\}$\textup{)}:
$$
	\rho(x_1)=-\la x_1,\quad \rho(x_2)=\la x_2, \quad \rho(y_1)=-2\la y_1,\quad \rho(y_2)=2\la y_2, \quad \text{where $\lambda\in \mathbb{K}^*$.}
	$$
\item[$(ii)$]  $(\g^4, \om, \rho)$ \textup{(}in the basis $\{x_1, x_2\mid  y_1, y_2\}$\textup{)}:
$$
	\rho(x_1)=\la x_1+\beta x_2,\quad \rho(x_2)=\mu x_1-2\la x_2, \quad \rho(y_1)=2\la y_1+\frac{\beta}{2} y_2,\quad \rho(y_2)=2\mu y_1-\lambda y_2,
	$$
 where $\lambda,\beta, \mu\in \mathbb{K}$ and $2\lambda^2+\beta\mu\neq 0$.
\end{itemize}
\end{Proposition}

\begin{proof}
Any flat quasi-Frobenius Lie superalgebra is isomorphic to one of the Lie superalgebras classified in Prop.~\ref{cls-dim4}. The other two Lie superalgebras, namely the first and the third, do not admit such a structure. Indeed, for these Lie superalgebras we have
$ [\g,\g]^\perp \neq Z(\g),$ and the claim follows from Prop.~\ref{roh-inv}.

Now we show that the second Lie superalgebra admits an even invertible $\omega$–antisymmetric endomorphism $\rho$ for which  $\rho\circ \ad_u$ is $\omega$-symmetric. We choose an even endomorphism defined on $\g$ by
$$
\rho(x_2)= a x_1 + b x_2,\qquad
\rho(y_2)= a_1 x_1 + b_1 x_2, \quad \text{where $a,b,a_1,b_1\in \mathbb{K}$.}
$$

$$
\omega(\rho [y_1,y_1],x_2)
= \omega(\rho(x_2),x_2)
= 2a
= -\omega(y_1,\rho[y_1,x_2])
= 0,
$$
hence $a=0$. Similarly,
$$
\omega(\rho [y_1,y_1],x_1)
= \omega(\rho(x_2),x_1)
= -2b
= -\omega(y_1,\rho[y_1,x_1])
= \omega(y_1,\rho(y_2))
= b_1,
$$
which gives $b_1 = -2b$. Next,
$$
\omega(\rho[x_1,y_1],x_1)
= \omega(\rho(y_2),x_1)
= 2a_1
= \omega(y_1,\rho[x_1,x_1])
= 0,
$$
thus $a_1=0$. Therefore,
$$
\rho(x_2)= b x_2,\qquad
\rho(y_2)= -2b y_2.
$$
Since $\rho$ is invertible, then $b\neq0$. As $\rho$ is $\omega$-antisymmetric, we deduce
$$
\rho(x_1)= -b x_1,\qquad
\rho(y_1)= 2b y_1.
$$
We show in the same way that the Lie superalgebra $(\g^4,\omega)$ also admits  an even invertible $\omega$-antisymmetric endomorphism $\rho$ satisfying  $\rho\circ \ad_u$ is $\omega$-symmetric.

Finally, when $\rho$ is odd, no such structure exists for the remaining nonabelian Lie superalgebra, according to Prop.~\ref{rho-odd}.
This completes the proof.
\end{proof}

\ssbegin{Example}
Let $(\mathfrak{b}:=\mathfrak{b}_{\bar{0}}, \om_\mathfrak{b})$ be a $4$-dimensional abelian Lie algebra equipped with a quasi-Frobenius  form $\omega_{\mathfrak{b}}$ given in the basis $\{e_1,e_2,e_3,e_4\}$ by 
$$
\omega_{\mathfrak{b}}=e_1^*\wedge e_4^*+e_2^*\wedge e_3^*.
$$
Consider the even linear map $\rho_{\mathfrak{b}}$ defined by
$$
\rho_{\mathfrak{b}}(e_1)=-2\lambda e_1,\;\rho_{\mathfrak{b}}(e_2)=2\lambda e_2,\;
\rho_{\mathfrak{b}}(e_3)=-2\lambda e_3,\;\rho_{\mathfrak{b}}(e_4)=2\lambda e_4, \quad \text{where $\lambda \in \mathbb{K}\backslash \{0\}$.}
$$
It is clear that $\rho_{\mathfrak{b}}$ is invertible and $\omega_{\mathfrak{b}}$-antisymmetric.

Choose $b_0=0$ and $a_0=0$. Let $\xi$ be an even endomorphism of $\mathfrak{b}$ defined by
\[
	\xi=a E_{1,3}-a E_{2,4}, \quad \text{where $a\neq 0$.}
\]
It is easy to check that $(\xi, b_0, a_0, \lambda)$  satisfies Systems~\eqref{eq:claim2} and~\eqref{DEX1}. Let $V:=\mathbb{K} d$  be a
one-dimensional vector space and $V^* = \mathbb{K} e$ be its dual. Then $\mathfrak{g}=V \oplus \mathfrak{b}\oplus V^*$ is an even flat double extension of $(\mathfrak{b}, \omega_{\mathfrak{b}}, \rho_{\mathfrak{b}})$ by means of  $(\xi, b_0, a_0, \lambda)$. According to Theorem \ref{THQ1}, the Lie bracket, the  form $\om$, and the  endomorphism $\rho$ on $\mathfrak{g}$ are defined as follows: (in  the basis $\{d, e , e_1, e_2, e_3, e_4\}$)
\begin{itemize}
\item[$\bullet$] Lie bracket is given by  
$$[d, e_3]=a e_1,\quad [d, e_4]=-a e_2,\quad [e_3, e_4]=2a e,$$
\item[$\bullet$] The form $\om$ is given by 
$$\om=e^*\wedge d^*+e_1^*\wedge e_4^*+e_2^*\wedge e_3^*,$$
\item[$\bullet$] The endomorphism $\rho$ is given by 
  $$
  \rho(e) = \lambda e, \; \rho(d) = -\lambda d, \;
  \rho(e_1) = -2 \lambda e_1, \; \rho(e_2) = 2 \lambda e_2, \;
  \rho(e_3) = -2 \lambda e_3, \; \rho(e_4) = 2 \lambda e_4.
  $$

\end{itemize}
With these definitions, $(\mathfrak{g}, \omega, \rho)$ is a flat QQF-algebra.

\end{Example}

\ssbegin{Example}
Let $(\mathfrak{b}:=\mathfrak{b}_{\bar{1}}, \om_\mathfrak{b})$ be a $4$-dimensional purely odd abelian Lie superalgebra equipped with an  orthosymplectic form $\omega_{\mathfrak{b}}$ given in the basis $\{0\, |\,  e_1,e_2,e_3,e_4\}$ by
$$
\omega_{\mathfrak{b}}= -e_1^*\wedge e_4^*-e_2^*\wedge e_3^*.
$$
Consider the even linear map $\rho_{\mathfrak{b}}$ defined in this basis by
$$
\rho_{\mathfrak{b}}(e_1)=-2\lambda e_1,\;\rho_{\mathfrak{b}}(e_2)=-2\lambda e_2,\;
\rho_{\mathfrak{b}}(e_3)=2\lambda e_3,\;\rho_{\mathfrak{b}}(e_4)=2\lambda e_4, \quad \text{where $\lambda \in \mathbb{K}^*$.}
$$
It is clear that $\rho_{\mathfrak{b}}$ is invertible and $\omega_{\mathfrak{b}}$-antisymmetric.

Choose $b_0=a_0=0$. Let $\xi$ be an even endomorphism of $\mathfrak{b}$ defined by 
\[
	\xi=a E_{1,3}+ a E_{2,4}, \quad \text{where $a\neq 0$.}
\]
It is easy to check that $(\xi, b_0, a_0, \lambda)$  satisfies Systems~\eqref{eq:claim2} and~\eqref{DEX1}. Let $V:=\mathbb{K} d$  be a
one-dimensional vector space and $V^* = \mathbb{K} e$ be its dual. Then $\mathfrak{g}=V \oplus \mathfrak{b}\oplus V^*$ is even flat double extension of $(\mathfrak{b}, \omega_{\mathfrak{b}}, \rho_{\mathfrak{b}})$ by means of  $(\xi, b_0, a_0, \lambda)$. According to Theorem \ref{THQ1}, the Lie bracket, the even form $\om$, and the even endomorphism $\rho$ on $\mathfrak{g}$ are defined as follows:  (in the basis $\{d, e \mid e_1, e_2, e_3, e_4\}$)
\begin{itemize}
\item[$\bullet$] The Lie bracket is given by 
$$[d, e_3]=a e_1,\quad [d, e_4]=a e_2,\quad [e_3, e_4]=2a e,$$
\item[$\bullet$] The form $\om$ is given by 
$$\om=e^*\wedge d^*-e_1^*\wedge e_4^*-e_2^*\wedge e_3^*,$$
\item[$\bullet$] The even endomorphism $\rho$ is given by 
  $$
  \rho(e) = \lambda e, \quad \rho(d) = -\lambda d, \; 
  \rho(e_1) = -2 \lambda e_1, \;  \rho(e_2) = -2 \lambda e_2, \; 
  \rho(e_3) = 2 \lambda e_3, \;  \rho(e_4) = 2 \lambda e_4.
  $$

\end{itemize}
With these definitions, $(\mathfrak{g}, \omega, \rho)$ is a flat QQF-superalgebra where  $\rho$ is even.
\end{Example}

\ssbegin{Example}
Let $(\mathfrak{b}:=\mathfrak{b}_{\bar{0}}\oplus\mathfrak{b}_{\bar{1}}, \om_\mathfrak{b})$ be a $4$-dimensional  abelian Lie superalgebra equipped with an orthosymplectic form $\omega_{\mathfrak{b}}$ given in the basis $\{e_1,e_2\mid e_3,e_4\}$ as follows: 
$$
\omega_{\mathfrak{b}}=e_1^*\wedge e_2^*- e_3^*\wedge e_4^*.
$$
Consider the even linear map $\rho_{\mathfrak{b}}$ defined by
$$
\rho_{\mathfrak{b}}(e_1)=\lambda e_1,\;\rho_{\mathfrak{b}}(e_2)=-\lambda e_2,\;
\rho_{\mathfrak{b}}(e_3)=-2\lambda e_3,\;\rho_{\mathfrak{b}}(e_4)=2\lambda e_4, \quad \text{where $\lambda \in \mathbb{K}\backslash \{0\}$.}
$$
It is clear that $\rho_{\mathfrak{b}}$ is invertible and $\omega_{\mathfrak{b}}$-antisymmetric.

Choose $b_0=a_0=0$. Let $\xi$ be an even endomorphism of $\mathfrak{b}$ defined by
\[
	\xi=a E_{3,4}, \quad \text{where $a\neq 0$.}
\]
It is easy to check that $(\xi, b_0, a_0, \lambda)$  satisfies Systems~\eqref{eq:claim2} and~\eqref{DEX1}. Let $V:=\mathbb{K} d$  be a
one-dimensional vector space and $V^* = \mathbb{K} e$ be its dual. Then $\mathfrak{g}=V \oplus \mathfrak{b}\oplus V^*$ is even flat double extension of $(\mathfrak{b}, \omega_{\mathfrak{b}}, \rho_{\mathfrak{b}})$ by means of  $(\xi, b_0, a_0, \lambda)$. According to Theorem \ref{THQ1}, the Lie bracket, the even form $\om$, and the even endomorphism $\rho$ on $\mathfrak{g}$ are defined as follows: (in the basis $\{d, e, e_1, e_2\mid e_3, e_4\}$)
\begin{itemize}
\item[$\bullet$] The Lie bracket is given by  
$$[d, e_4]=a e_3,\quad [e_4, e_4]=2a e,$$
\item[$\bullet$] The form $\om$ is given by 
$$\om=e^*\wedge d^*+e_1^*\wedge e_2^*-e_3^*\wedge e_4^*,$$
\item[$\bullet$] The even endomorphism $\rho$ is given by
  $$
  \rho(e) = \lambda e, \quad \rho(d) = -\lambda d, \quad
  \rho(e_1) = \lambda e_1, \quad \rho(e_2) = - \lambda e_2, \quad
  \rho(e_3) = -2 \lambda e_3, \quad \rho(e_4) = 2 \lambda e_4.
  $$

\end{itemize}
With these definitions, $(\mathfrak{g}, \omega, \rho)$ is a flat QQF-superalgebra where  $\rho$ is even.
\end{Example}

\ssbegin{Example}
Let $(\mathfrak{b}:=\mathfrak{b}_{\bar{0}}\oplus\mathfrak{b}_{\bar{1}}, \om_\mathfrak{b})$ be a $4$-dimensional  abelian Lie superalgebra equipped with a periplectic form $\omega_{\mathfrak{b}}$ given in the basis $\{e_1,e_2\mid e_3,e_4\}$ by
$$
\omega_{\mathfrak{b}}=e_1^*\wedge e_3^*+e_2^*\wedge e_4^*.
$$
Consider the even linear map $\rho_{\mathfrak{b}}$ defined by
$$
\rho_{\mathfrak{b}}(e_1)=2\la e_1,\;\rho_{\mathfrak{b}}(e_2)=-2\lambda e_2,\;
\rho_{\mathfrak{b}}(e_3)=-2\la e_3,\;\rho_{\mathfrak{b}}(e_4)=2\la e_4, \quad \text{where $\lambda \in \mathbb{K}\backslash \{0\}$.}
$$
It is clear that $\rho_{\mathfrak{b}}$ is invertible and $\omega_{\mathfrak{b}}$-antisymmetric.

Choose $b_0=a_0=0$. Let $\xi$ be an even endomorphism of $\mathfrak{b}$ defined by
\[
	\xi=a E_{2,1}+a E_{3,4}, \quad \text{where $a\neq 0$.}
\]
It is easy to check that $(\xi, b_0, a_0, \lambda)$  satisfies Systems~\eqref{eq:claim2} and~\eqref{DEX1}. Let $V:=\mathbb{K} d$  be a
one-dimensional vector space and $V^* = \mathbb{K} e$ be its dual. Then $\mathfrak{g}=V \oplus \mathfrak{b}\oplus V^*$ is even flat double extension of $(\mathfrak{b}, \omega_{\mathfrak{b}}, \rho_{\mathfrak{b}})$ by means of  $(\xi, b_0, a_0, \lambda)$. According to Theorem \ref{THQ1}, the Lie bracket, the odd form $\om$, and the even endomorphism $\rho$ on $\mathfrak{g}$ are defined as follows: (in the basis $\{d, e, e_1, e_2\mid e_3, e_4\}$)
\begin{itemize}
\item[$\bullet$] The Lie bracket is given by 
$$[d, e_1]=-a e_2,\quad [d, e_4]=-a e_3,\quad [e_1, e_4]=2a e,$$
\item[$\bullet$] The form $\om$ is given by 
$$\om=e^*\wedge d^*+e_1^*\wedge e_3^*+e_2^*\wedge e_4^*,$$
\item[$\bullet$] The even endomorphism $\rho$ is given by 
  $$
  \rho(e) = \lambda e, \quad \rho(d) = -\lambda d, \quad
  \rho(e_1)=2\la e_1,\;\rho(e_2)=-2\lambda e_2,\;
\rho(e_3)=-2\la e_3,\;\rho(e_4)=2\la e_4.
  $$
\end{itemize}
With these definitions, $(\mathfrak{g}, \omega, \rho)$ is a flat QQF-superalgebra where  $\rho$ is even.
\end{Example}

We established a classification of flat quasi-quadratic quasi-Frobenius QQF-superalgebras of dimension four when $\rho$ is even, and we provided several low-dimensional examples, in particular in dimension six.
Moreover, Corollaries \ref{cosuite} and \ref{cosuite1} show that every flat QQF-superalgebra for which $\rho$ is odd necessarily has dimension $4n$. In addition, Prop. \ref{rho-odd} proves that any 4-dimensional flat QQF-superalgebra for which $\rho$ is odd must be abelian.
Consequently, the smallest dimension in which non-abelian examples may occur in the odd case is eight. In order to illustrate this situation, and motivated by the theory of double extensions, we now construct an explicit example of a flat QQF-superalgebra of dimension eight obtained via a suitable double extension.

\ssbegin{Example}
Let $(\mathfrak{b}:=\mathfrak{b}_{\bar{0}}\oplus\mathfrak{b}_{\bar{1}}, \om_\mathfrak{b})$ be a $4$-dimensional  abelian Lie superalgebra equipped with a  periplectic form $\omega_{\mathfrak{b}}$ given in the basis $\{f_1,f_2\mid f_3,f_4\}$ by
$$
\omega_{\mathfrak{b}}=f_1^*\wedge f_3^*+f_2^*\wedge f_4^*.
$$
Consider the odd linear map $\rho_{\mathfrak{b}}$ defined by
$$
\rho_{\mathfrak{b}}(f_1)=-2\la f_4,\;\rho_{\mathfrak{b}}(f_2)=-2\lambda f_3,\;
\rho_{\mathfrak{b}}(f_3)=\la f_2,\;\rho_{\mathfrak{b}}(f_4)=-\la f_1, \quad \text{where $\lambda \in \mathbb{K}\backslash \{0\}$. }
$$
It is clear that $\rho_{\mathfrak{b}}$ is invertible and $\omega_{\mathfrak{b}}$-antisymmetric.

Choose $b_0=b_1=c_0=c_1=a_0=a_1=0$. Let $\xi_0$ be an even endomorphism of $\mathfrak{b}$, and let $\xi_1$ be an odd endomorphism of $\mathfrak{b}$, defined by 
\[
	\xi_0=2a E_{2,1}+a E_{3,4}, \quad \xi_1=\frac{3}{2}aE_{2,4}-a E_{3,1}, \quad \text{where $a\neq 0$.} 
\]
It is easy to check that $(\xi_0,\xi_1, 0,0,0,0,0,0)$ satisfying the system \eqref{2plane1}-- \eqref{2plane6} and \eqref{sy-rho1}. Let $V:=\mathbb{K} d_0\oplus \mathbb{K} d_1$  be a 
two-dimensional vector space and $\Pi(V^*) = \mathbb{K} e_0\oplus \mathbb{K} e_1$. Then $\mathfrak{g}=V \oplus \mathfrak{b}\oplus \Pi(V^*)$ is a flat planar double extension of $(\mathfrak{b}, \omega_{\mathfrak{b}}, \rho_{\mathfrak{b}})$ by means of  $(\xi_0,\xi_1, 0,0,0,0,0,0)$. According to Theorems \ref{THDOUBLE} and \ref{THRE2}, the Lie bracket, the odd form $\om$, and the odd endomorphism $\rho$ on $\mathfrak{g}$ are defined as follows (in the basis $\{d_0, e_0, f_1, f_2\mid d_1, e_1, f_3, f_4\}$):
\begin{itemize}
\item[$\bullet$] The Lie bracket is given by 
$$[d_0, f_1]=-3 a f_2,\quad [d_1, f_1]=3 a f_3,\quad [d_1, f_4]=\frac{3}{2}a f_2,\quad [f_1, f_4]=3a e_1, \quad [f_4, f_4]=-3 a e_0,$$
\item[$\bullet$] The form $\om$ is given by 
$$\om=e_0^*\wedge d^*_1+e_1^*\wedge d^*_0+f_1^*\wedge f_3^*+f_2^*\wedge f_4^*,$$
\item[$\bullet$] The odd endomorphism $\rho$ is given by 
  \begin{align*}
&\rho(e_0)=\la e_1, \quad \rho(e_1)=\la e_0,\quad \rho(d_0)=-\la d_1,\quad \rho(d_1)=\la d_0,\quad  \rho(f_1)=-2\la f_4,\\&\rho(f_2)=-2\lambda f_3,\quad
\rho(f_3)=\la f_2,\;\rho(f_4)=-\la f_1,
 \end{align*}
\end{itemize}
With these definitions, $(\mathfrak{g}, \omega, \rho)$ is a flat QQF-superalgebra where  $\rho$ is odd.
\end{Example}

{\bf Acknowledgments:} We thank the referee for helping us  improve the presentation of the paper. 


\end{document}